\documentclass[11pt]{amsart}
\usepackage{amssymb}
\usepackage{amsmath}
\usepackage{amscd}
\usepackage{amsthm}
\usepackage{epsfig}  
\usepackage{epsfig}  
\usepackage{epsfig}  
\usepackage{subfigure}
\usepackage{changebar}
\usepackage{wrapfig,mathrsfs}
\usepackage{euscript,euler, latexsym, amssymb}
\usepackage{delarray}
\usepackage{commands}
\usepackage[all]{xy}
\numberwithin{equation}{section}
\begin{document}  

\title[The unstable integral homology of the mapping class groups]{The unstable integral homology of the mapping class groups of 
a surface with boundary.}
\author{V\'eronique Godin}
\address{Department of Mathematics\\
  Harvard University\\
   Cambridge, MA 02138}
\email{godin@math.harvard.edu}
\date{\today}

\begin{abstract}
We construct a graph complex calculating the \textbf{integral} homology of the bordered mapping class groups. We compute the homology of the bordered mapping class groups of the surfaces $S_{1,1}$, $S_{1,2}$ and $S_{2,1}$.  Using the circle action on this graph complex, we build a double complex and a spectral sequence converging to the homology of the unbordered mapping class groups. We compute the homology of the punctured mapping class groups associated to the surfaces $S_{1,1}$ and $S_{2,1}$. Finally, we use Miller's operad to get the first Kudo-Araki and Browder operations on our graph complex. We also consider an unstable version of the higher Kudo-Araki-Dyer-Lashoff operations.
\end{abstract}

\maketitle

\newcounter{temp}

\section{Introduction}
Let $S=S_{g,n}$ be a surface of genus $g$ with $n$ boundary components. We assume that the boundary of $S$ is not empty. The \emph{bordered mapping class group}  of $S$
\[ \cM(S;\p) =\pi_0 \Diff(S;\p S)\] 
is the group of isotopy classes of \mbox{self-}diffeo\-mor\-phisms which
fix the boundary pointwise. The \emph{punctured mapping class group} of $S$ 
\[\cM(S) = \pi_0\Diff(S,\p_1S,\dots\p_nS)\]
is the group of isotopy classes of orientation-preserving diffeomorphisms which restrict to diffeomorphisms of each boundary components. 

These mapping class groups act on the appropriate Teichm\"uller spaces with quotient the moduli space of complex surfaces of the suitable type. The action of $\cM(S;\p)$ is free and the integral homology of the group $\cM(S;\p)$ is the integral homology of the moduli space of conformal structure on $S$ with cylindrical ends. Since the action of $\cM(S_{g,n})$ has only finite isotropy, the rational homology of the group $\cM(S_{g,n})$ coincide with the rational homology of the moduli space of conformal structure on a closed surface $S_{g,0}$ with $n$ marked points.

By gluing a twice-punctured torus to the unique boundary of a surface $S_{g,1}$ with genus $g$, we get group homomorphisms
\[\cdots\map\cM(S_{g,1};\p)\maplu{\psi_g} \cM(S_{g+1,1};\p)\maplu{\psi_{g+1}} \cM(S_{g+2,1};\p)\map\cdots\]
whose direct limit we denote by $\cM_{\infty}$. Harer \cite{Harer_Stability} showed that the $\psi_g$'s induce cohomological and homological isomorphisms in a range of dimensions increasing with the genus $g$. In this stable range, the homology and the cohomology groups
\[H_k(\cM(S_{g,1};\p)) \cong H_k(\cM_{\infty})\qquad H^k(\cM(S_{g,1};\p))\cong H^k(\cM_\infty) \qquad g\geq 2k+1\]
are independent of the genus and are called the stable homology and cohomology of the bordered mapping class group. Mumford conjectured in \cite{Mumford} that the stable rational cohomology of these mapping class groups is a polynomial algebra 
\[ H^*(\cM_{\infty};\bQ) \cong \bQ[\kappa_1,\kappa_2,\ldots].\] 
in the tautological classes which are obtained from the Chern class of the vertical tangent bundles of the universal surface bundle. 

By gluing surfaces $S_{g_i,1}$ to the first $n$ boundary components of the generalized pair of pants $S_{0,n+1}$, Miller defined homomorphisms
\begin{equation}\label{eq:framed}
\cM(S_{g_1,1};\p) \times \cdot \cM(S_{g_n,1};\p)\times \cM(S_{0,n};\p)\map \cM(S_{\sum{g_i},1};\p).
\end{equation}
Using a recognition principle, he showed that the group completion of
\[X=\Du B\cM(S_{g,1};\p)\]
is a two-fold loop space.
In \cite{Tillmann}, Tillmann used a cobordism category to extend Miller's result. She showed that the group completion of $X$ has the cohomology of an infinite loop space. Madsen and Weiss \cite{Mad} then identified this infinite loop space to be $\Omega^{\infty}\cC P^{\infty}_{-1}$.  Rationally, this result proves the Mumford conjecture. However, it also gives information about the torsion of the stable cohomology. Using the Madsen-Weiss theorem, Galatius in \cite{Galatius_Modp} computed the mod-$p$ cohomology of the infinite loop space $\Omega^{\infty}\cC P^{\infty}_{-1}$, uncovering a rich and unexpected torsion component of the stable cohomology of the mapping class groups.  

These recent exciting results describe the stable cohomology of the mapping class groups. However, not much is known about their unstable cohomology. The goal of this paper is to study this integral unstable homology.

Although the unstable homology is interesting in itself, we have a specific application in mind. In \cite{Cohen_Godin}, Cohen and the author have used a special type of fat graphs to define operations on the free loop space $LM$ of an orientable manifold extending earlier work of Chas and Sullivan \cite{Chas_Sullivan}. It is conjectured that these operations extend to operations 
\[H_*(\cM(S_{p+q};\p)) \tens H_*(LM)^{\tens p} \map H_*(LM)^{\tens q}\]
parameterized by the homology of the bordered mapping class groups. In particular, this would give that the Chas and Sullivan product is part of an $E_\infty$-structure on $H_*(LM)$. Antonio Ram\'irez and the author are working on using the model introduced in this paper to define these operations. Also the homology classes found in this paper will give example of these operations on $H_*(LM)$ and may prove essential to show the non-triviality of these higher homological operations.

\subsection{Fat graphs and the mapping class groups.}

A fat graph or ribbon graph is a finite connected graph with a cyclic ordering of the half-edges incident to each vertex. From a fat graph, we construct a surface by replacing each edge by a thin ribbon and by gluing these ribbons at the vertices according to the cyclic orderings.

Following ideas of Thurston, Strebel \cite{Strebel}, Bowditch and Epstein \cite{Bowditch} and Penner \cite{Penner} constructed a triangulation of the decorated Teichm\"uller space of a punctured Riemann surface $S$ which is equivariant under the action of the unbordered mapping class group. 
The quotient space, in which a point is an isomorphism class of metric fat graphs, gives a model for the corresponding decorated moduli space.   

The spaces of fat graphs are filtered by the combinatorics of the graphs. This stratification has been used by Penner \cite{Penner}, Harer and Zagier \cite{Harer_Zagier} to compute the Euler characteristic of the moduli space and by Kontsevich to prove Witten's  conjecture about the intersection numbers of the Miller-Morita-Mumford classes in the Deligne-Mumford compactification of the moduli space of punctured surfaces \cite{Kontsevich_Intersection}.

The spectral sequence associated to the combinatorial filtration collapses rationally to a complex. This complex is generated freely by isomorphism classes of oriented fat graphs and its boundary maps can be described combinatorially. Penner \cite{Penner_complex} was first to build this graph complex $\cG\cC^*$ which by construction computes the rational cohomology of the moduli space of marked Riemann surfaces. 

\subsection{Results.}
To study bordered mapping class group, we first extend the notion of fat graph.
In section \ref{sec:proof of thm}, we define a bordered fat graph to be a fat graph with exactly one leaf (vertex with a single edge attached
to it) for each boundary component. Each of these leaves gives a marked point on the corresponding boundary
component. We denote by $\Fatb$ the category of isomorphism classes of bordered fat graphs. Using the work of Harer \cite{Harer}, we show that there is an homotopy equivalence
\[(|\Fatb|)_+ \simeq \Du_{g\geq0 ,n\geq 1} B\cM(S_{g,n};\p)\]
between the geometric realization of $\Fatb$ with an added base point and the classifying spaces of the bordered mapping class groups.

Following Penner, we define, in section \ref{sec:CW}, a combinatorial filtration on our categorical model $\Fatb$.
We show in theorem \ref{thm:fp/fp-1} that this filtration is the skeleton of a CW-structure on $|\Fatb|$ with exactly one cell
for each isomorphism class of bordered fat graphs.  
The cellular chain complex $\GC_*$ of this CW-structure gives the equivalent of the graph complex for
bordered fat graphs. Note that although the original graph complex gives an isomorphism only on the rational
homology, we get an integral result
\[H_*(\GC_*;\bZ)\cong\bigoplus\ H_*(\cM(S;\p);\bZ).\]
Using $\GC_*$, we compute the homology of $\cM(S_{1,1};\p)$ directly.
\[H_*(\cM(S_{1,1};\p))\cong \begin{cases}
	\bZ&*=0,1\\
	0 &*\geq 2
	\end{cases}.\]
Using a computer algebra, we get
\[H_*(\cM(S_{1,2};\p))\cong\begin{cases}
\bZ &*=0\\
\bZ\oplus\bZ & *=1\\
\bZ/2\oplus \bZ & *=2\\
\bZ/2 & *=3\\
0 & *\geq 4
\end{cases}\quad
H_*(\cM(S_{2,1};\p)) \cong\begin{cases}
\bZ & *=0\\
\bZ/10& *=1\\
\bZ/2 & *=2\\
\bZ/2\oplus \bZ & *=3\\
\bZ/6 &*=4\\
0 &*\geq 5.
\end{cases}\]
The results for $\cM(S_{1,1};\p)$ and $\cM(S_{2,1};\p)$ match the computations of Ehrenfried in \cite{ehr}.

To get at the integral homology of the punctured mapping class group of a surface $S=S_{g,n}$, we consider the exact sequence of topological groups
\[0\map \Diff(S;\p) \map \Diff(S,\p_1S, \dots,\p_nS) \map \prod \Diff(\p_i S,\p_i S)\map 0\]
whose homotopy long exact sequence gives
\[0\map \bZ^n\map \cM(S;\p)\map \cM(S) \map 0.\]
In fact, the group $\cM(S)$ is obtained from $\cM(S;\p)$ by killing the Dehn twists around the boundary components.  Since these Dehn twists are central, the bordered mapping class group is a central extension of the punctured one. 

For a surface with a single boundary component, this extension gives a double complex structure on the vector space $\GC_*\tens \bZ[u]$ which calculates the integral homology of the unbordered mapping class groups $\cM(S_{g,n})$. Using this double complex and its associated spectral sequence, we compute the following homology groups.
\begin{eqnarray*}
H_*(\cM(S_{1,1});\bZ) &\cong& \begin{cases}
\bZ&*=0\\
\bZ/12 &*=2k+1\\
0 &*=2k+2.
\end{cases}\\
H_*(\cM(S_{2,1});\bZ) &\cong& \begin{cases}
 \bZ &*=0\\
\bZ/10 &*=1\\
\bZ/2\oplus \bZ&*=2\\
\bZ/2\oplus \bZ/120 \oplus \bZ/10&*=2k+3\\
\bZ/2\oplus\bZ/6 &*=2k+4
\end{cases}
\end{eqnarray*}

We then translate Miller's homomorphism to our models. We get a product on the bordered graph complex
\[\GC_{p}\tens\GC_{q}\map\GC_{p+q}.\]
We know that this product is homotopy commutative. Using this homotopy, we define maps
\[\GC_{p}\map \GC_{2p+1}\qquad \GC_{p}\tens \GC_{q}\maplu{\phi} \GC_{p+q+1}\]
which induce the first Araki-Kudo and Browder operations at the homology level. The Browder operations are obstruction to an $n$-loop space being an higher loop space. Since the infinite loop space of Tillmann extends Miller's double loop space, the Browder operation hit only unstable classes. The result of Tillmann also gives higher Araki-Kudo-Dyer-Lashoff operations                                                                                                                                                                                        
\[Q_{i,p}:H_k(\cM_\infty;\bZ/p)\map H_{pk+i}(\cM_\infty;\bZ/p).\]
Using an idea of Cohen and Tillmann \cite{Cohen_Tillmann}, we build operations
\[ \widetilde{Q}_{i,p}: H_k(\cM(S_{g,1};\p);\bZ/p) \map H_{pk+i}(\cM(S_{pg,1};\bZ/p).\]
The $Q_{i,p}$ are obtained from projecting to the stable bordered mapping class group whose homology is a direct summand of the stable homology of the unbordered mapping class group.

\subsection{Remark.}
There exists other models for the classifying spaces of the bordered mapping class groups and some computations in low genus have already been made. In \cite{Harer} Harer extended the notion of arc complexes of Strebel to bordered surfaces. This model was subsequently used by Kaufmann, Livernet and Penner \cite{Kauf_Liv_Pen} to define an operad structure on a compactification of the moduli space of bordered Riemann surfaces. Our proof that $\Fatb$ realizes to a classifying space for the bordered mapping class group relies heavily on the work of Harer.

In \cite{Bodig} B\"odigheimer has constructed a configuration space model $\mathfrak{Rad}$ consisting of pairs of radial slits on annuli. Using this model, Ehrenfried in \cite{ehr} calculated the integral homology of $H_*(\cM(S_{1,1};\p))$ and $H_*(\cM(S_{2,1};\p)$. His computations agree with ours.

Although the model introduced in this paper is close to the arc complex model, both its categorical nature and its use of fat graphs will prove essential for future applications to string topology. A category can be studied by using techniques of algebraic topology and homotopy theory and future applications will utilize homotopy limits, techniques of algebraic $K$-theory, symmetric monoidal categories and infinite loop spaces, as well as theorems of McDuff, Segal, Quillen and Grothendiek.

The author would like to thank Ralph Cohen, Daniel Ford, Tyler Lawson, Antonio Ram\'irez and Ralph Kaufmann for interesting conversations on these topics.

\section[A categorical model $\Fatb$ for the bordered mapping class group.]{A categorical model for the bordered mapping class group.}

Let $S=S_{g,n}$ be a surface of genus $g$ with $n$ boundary components. Again 
\[\cM(S;\p) =\pi_0 \Diff(S;\p)\]
will be the group of isotopy classes of orientation-preserving diffeomorphisms of $S$ that fixes the
boundary pointwise. In this section, we construct a category $\Fatb$ whose objects are isomorphism
classes of bordered fat graphs. The geometric realization of this category
is a classifying space for the bordered mapping class groups. More precisely,
\[|\Fatb|\simeq \Du_{g\geq0, n\geq1} B\cM(S_{g,n};\p) \quad (g,n)\neq (0,1).\]

\subsection{Fat graphs and punctured Riemann surfaces.}

We briefly introduce the classical fat-graph model for the punctured mapping class groups.

\begin{definition}
A \textbf{combinatorial graph} $G=(V,H,s,i)$ consists of a set of vertices $V$, a set of half-edges $H$, a map $s:H\smap V$ and an involution
$i:H\smap H$ without fixed points.
\end{definition}

The map $s$ sends an half-edge to its source. 
The involution $i$ pairs an half-edge with its other half and an edge of $G$ is an orbit of the involution $i$.
The geometric realization of $G$ is a CW-complex $|G|$ with vertices $V$, 1-cells $|G|_1=H/i$
and no higher cells.

A tree $T=(V_T,H_T)$ of $G$ is a subgraph of $G$ whose geometric realization $|T|$ is contractible.
A forest $F$ of $G$ is a subset of the set of half-edges of $G$ that is closed under $i_G$ such that the realization
$|F|$ of the combinatorial graph 
\[(V_G, F, i|_{F}, s|_{F})\]
is a disjoint union of contractible spaces. 

We extend the maps $i$ and $s$ to be the identity on vertices and define
a \textbf{morphism of combinatorial graph} $\phi : G \maps \tildeG$ to be a map of sets 
\[\phi : V_G\du\, H_G\maps V_{\tildeG} \du\, H_\tildeG\] 
that commutes with $s$ and $i$ and such that 
\begin{enumerate}
\item $\phi^{-1}(v)$ is a tree for every vertex $v \in V_\tildeG$
\item $\phi^{-1}(A)$ contains a single half-edge of $G$ for every half-edge $A$ of $\tildeG$
\end{enumerate}
Such a morphism of combinatorial graph induces a simplicial and surjective homotopy equivalence on the geometric
realizations.

\begin{definition}
A \textbf{fat graph} $\G$ is a combinatorial graph $G$ together with a cyclic ordering $\s_v$ of the half-edges incident to each
vertex $v$.
\end{definition}

Since each half-edge is incident to one vertex, the cycles $\s_v$ give a permutation $\s=(\s_v)$ on the set of half-edges, which we call the \emph{fat structure}. It sends an half-edge to its successor in the cyclic ordering at its source vertex.  

\begin{figure} 
\begin{center}
\mbox{
\subfigure[$\G$\label{fig:firstfatgraph}]{\epsfig{file=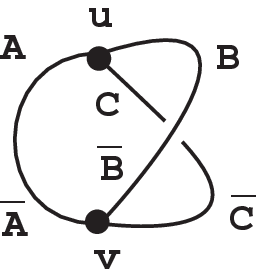, width=65pt}}\qquad
\subfigure[$\Sigma_{\G}$]{\epsfig{file=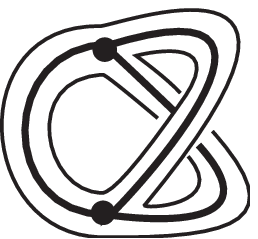, width=65pt}}\qquad\qquad
\subfigure[$\tilG$]{\epsfig{file=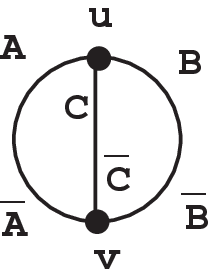, width=60pt}}\qquad
\subfigure[$\Sigma_{\tilG}$]{\epsfig{file=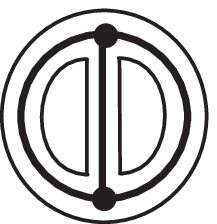, width=60pt}}}
\caption{Fat graphs and their associated surfaces}\label{fig:thickening}\label{fig:fat graph}
\end{center}
\end{figure}

\begin{example}
The combinatorial graph of figure \ref{fig:firstfatgraph} has a fat graph structure 
\[\s_u=(ABC) \quad \s_v=(\ol{ABC})\qquad \qquad \s=(ABC)(\ol{A}\ol{B}\ol{C}).\]
given by the clockwise orientation of the plane. 
\end{example}

A fat graph $\G=(G,\sigma)$ thickens to an oriented surface $\Sigma_\G$ with boundary. To build $\Sig_\G$, replace each
edge of $G$ by a strip and glue these strips together at the vertices according to the fat structure.

\begin{proposition}
The boundary components of $\Sig_\G$ correspond to the cycles of
$\om=\s\cdot i$
which is a permutation on the set of half-edges.
\end{proposition}

\begin{proof}
Any half-edge $A$ corresponds to one side of the strip of the edge $\{A,i(A)\}$ and hence corresponds to part of a boundary
component $\p_i\Sigma_\G$. Following this boundary along $A$ leads to its target vertex $v=s(i(A))$ where it follows the
next strip corresponding to $\sig(i(A))=\om(A)$ and so on.
\end{proof}

The orbits of the permutation $\om$ are called the \textbf{boundary cycles} of a fat graph. 
Since $\s=\om \cdot i$, the permutation $\om$ completely determines the fat structure.
\begin{example}
Changing the fat structure $\s$ affects the induced surface in a fundamental way.  In figure \ref{fig:fat graph}, the
fat graphs $\G$ and $\tilG$ have the same underlying graph but they have different fat structures. Their boundary cycles
are respectively,
\[\om_\G=(A\ol{B} C \ol{A} B\ol{C}) \qquad \om_\tilG=(A \ol{C})(B \ol{A})(C \ol{B}).\]
Hence the surface $\Sigma_{\tilG}$ has three boundary components and genus $0$ while the original $\Sigma_\G$ has a single boundary component
and genus $1$.
\end{example}

\label{sec:fat}
A morphism of fat graph 
\[\varphi : (G,\om)\maps(\tildeG,\tilom)\]
is a morphism of combinatorial graphs such that $\varphi(\om)=\tilom$. Hence $\tilom$ is obtained from $\om$ by replacing an half-edge
$A$ with $\varphi(A)$ if it is an half-edge or by skipping $A$ if $\varphi(A)$ is a vertex.
In \cite{Igu00}, Igusa built a category $\Fat$ whose objects are fat graphs with \textbf{no
univalent or bivalent vertices} and with ordered boundary cycles. The morphism of $\Fat$ are morphisms of fat graph which preserve the ordering of the boundary cyles. He
proved the following theorem by using ideas of Strebel \cite{Strebel}, Penner
\cite{Penner}, and Culler-Vogtman \cite{CV86}. 

\begin{theorem}[Igusa \cite{Igu00}]
The geometric realization of $\Fat$ is homotopy equivalent to the classifying space of the unbordered mapping class
groups 
\[|\Fat|\simeq\Du_{(g,n)\neq(0,1),(0,2)} B\cM(S_{g,n}).\]
Here $\cM(S_{g,n})$ is the group of isotopy classes of orientation-preserving self-diffeomorphisms of $S$ with no
boundary restrictions.
\end{theorem}

\subsection{A category of bordered fat graphs}
\begin{definition}\label{def:bivalent fat graphs}
A \textbf{bordered fat graph} is a fat graph $\Gb=(G,\om_1,\ldots\om_n)$ with ordered boundary components and \textbf{exactly} one leaf
(univalent vertex) in each boundary cycle. All other vertices of $\Gb$ are at least trivalent. 
\end{definition}

\begin{figure}
\begin{center}
\epsfig{file=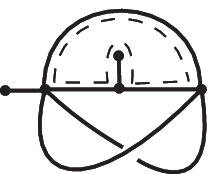, height=70pt}\caption{A bordered fat graph}\label{fig:bordered fat graph}
\end{center}
\end{figure}
\begin{example}
The bordered fat graph of figure \ref{fig:bordered fat graph} has two boundary components, one of which is dashed.
\end{example}

The edge leading to a leaf is called a leaf-edge. All other edges are interior edges.
The leaves determine a natural starting point for each boundary cycle. 
Note that the ordering of the boundary component gives an ordering of the leaves
and leaf-edges. 

A \textbf{morphism of bordered fat graphs}
\[\varphi : \Gb=(G,\om)\maps\tilGb= (\tildeG,\tilom)\] 
is a morphism of fat graphs sending the \ith leaf of $\G$ to the \ith leaf of $\tilG$.

We define $\Fatb$ to be the category  whose objects are isomorphism classes of bordered fat graphs. A
morphism $[\G^b_1]\smap [\G^b_2]$ of $\Fatb$ is an equivalence class of morphisms of bordered fat graphs where $\varphi\sim \tilvarphi$ if and
only if there are isomorphisms $\theta_i$ completing the following diagram.
\[\xymatrix{
\G^b_1 \ar[d]_{\theta_1}^{\cong}\ar[rr]^{\varphi} &&\G^b_2\ar[d]^{\theta_2}_{\cong}\\
\tilG^b_1 \ar[rr]^{\tilvarphi}&&\tilG^b_2
}\]

\begin{remark}\label{rmk:RFatb}
We could also define a category $\RFatb$ whose objects are
actual bordered fat graph (in some universe) and whose morphisms are
morphisms of bordered fat graph. It is not hard to show that the two
geometric realizations
\[ |\RFatb| \simeq |\Fatb| = Isom(\RFatb)\]
are homotopy equivalent.
\end{remark}

\begin{lemma}
$\Fatb$ is a well-defined category.
\end{lemma}
\begin{proof}
By definition an automorphism of bordered fat graph must preserve the first leaf
$L_1$. But then it must also preserve $A=\sigma(L_i)$ and so on. In fact, a bordered fat graph has
only the trivial automorphism and, similarly, two bordered fat graphs can have at most one
isomorphism between them.  

Two morphisms $[\psi_i:\G^b_i\smap\tilG^b_i]$ are composable if and only if $\tilG^b_1\cong\G^b_2$ and we define the
composition to be
\[\xymatrix{
\G^b_1\ar[r]^{\psi_1}\ar@/_1.1pc/[rrr]&\tilG^b_1\ar[r]_{\cong}^{\theta}&\G^b_2\ar[r]^{\psi_2}&\tilG^b_2
}\]
where $\theta$ is the unique isomorphism $\tilG^b_1\smap\G^b_2$.
\end{proof}

\begin{lemma}
Each morphism of $\Fatb$ is represented uniquely by a morphism
\[\psi_{[\Gb,F]} : [\Gb]\map[\Gb/F].\]
Here $[\Gb,F]$ is an isomorphism class of pairs $(\Gb,F)$ where $F$ is a forest
of $\Gb$ disjoint from the leaves. The bordered fat graph $\Gb/F$ is obtained
from $\Gb$ by collapsing each connected component of $F$ to a vertex 
\end{lemma}

\begin{proof}
For any pair $(\Gb,F)$, there exists a natural morphism of bordered fat graphs 
\[\phi_(\Gb,F): \Gb \smap \Gb/F\]
collapsing exactly the edges of $F$ which is the identity on the edges of $\Gb$
that are not in $F$.
Two such morphisms $\psi_{(\Gb,F)}$ and $\psi_{(\tilGb,\tilF)}$ are equivalent if and only if there exist isomorphisms completing
the following diagram.
\[\xymatrix{
\Gb\ar@{-->}[d]_{\theta}^{\cong}\ar[rr]_{\psi_{(\Gb,F)}}&&\Gb/F\ar@{-->}[d]^{\theta_F}_{\cong}\\
\tilGb\ar[rr]^{\psi_{(\tilGb,\tilF)}}&&\tilGb/\tilF
}\]
Any such $\theta$ defines an isomorphism of pairs $(\Gb,F)\smap(\tilGb,\tilF)$.

Given a morphism of bordered fat graphs $\varphi:\Gb\smap\tilGb$, let $F := \du \varphi^{-1}(v)$.
Define a morphism $\theta:\Gb/F\smap\tilGb$ by sending an edge $A$ not in $F$ to the edge $\varphi(A)$. Since the
diagram 
\[\xymatrix{
\Gb \ar[rr]^{\psi_{(\Gb,F)}} \ar[drr]_{\varphi}&&\Gb/F\ar[d]^{\theta}_\cong\\
&&\tilGb
}\]
commutes, $\varphi$ is equivalent to $\psi_{(\Gb,F)}$. 
\end{proof}

Denote by $\Fatbgn$ the full subcategory of $\Fatb$ which contains the
isomorphism classes of bordered fat graph $[\Gb]$ whose surface $\Sig_\Gb$ have
genus $g$ and $n$ boundary components. We will also call $g$ the genus of $\Gb$.
A morphism of $\Fatb$ induces an homotopy equivalence on the surface which preserve the number of boundary components and 
\[\Fatb=\Du \Fatbgn.\]

\subsection{Bordered fat graphs and the mapping class groups}
\label{sec:proof of thm}

\begin{theorem}\label{thm:M(S;p)}\label{thm:fatb}
Let $S_{g,n}$ be a surface of genus $g$ with $n$ boundary components. There is an homotopy equivalence
\[|\Fatb|\simeq \Du_{g\geq0 ,n\geq 1} B\cM(S_{g,n};\p)\quad (g,n)\neq(0,1)\]
between the geometric realization of $\Fatb$ and the classifying of the bordered mapping class groups.
\end{theorem}

For a fix $g$ and $n$, we now consider the category $E\Fatbgn$ of marked bordered fat graphs and its
relationship with the arc complex of the surface $S=S_{g,n}$. Fix an ordering of the boundary components of $S$ and fix a marked point
$x_i$ on each boundary component $\p_iS$. Given a bordered fat graph $\Gb$, let $v_i$ denote its \ith leaf.

\begin{definition}\label{def:marking}
A {\bf marking} $[H]$ of the bordered fat graph $\Gb$ is an isotopy class of embeddings $H:|\Gb|\smap S_{g,n}$ such that $H(v_i)=x_i$
and such that the cyclic ordering $\sigma_v$ at a vertex $v$ of $\Gb$ is induced by the orientation on $S$ at $H(v)$.
\end{definition}

Define $E\Fatbgn$ to be the category whose objects are isomorphisms classes of marked bordered fat graphs $[\Gb,[H]]$. $E\Fatbgn$ has one morphism
\[\psi_{[\Gb,F]}:[\Gb,[H]]\map[\Gb/F,[H_F]]\]
for each forest $F$ of $\Gb$ disjoint from the leaves.
Here $[H_F]$ is obtained from the marking $[H]$ on $\Gb$ by collapsing the edges of $F$ so that
\[\xymatrix{
|\Gb|\ar[rr]^-{|\psi_{[\Gb,F]}|}&&|\Gb/F|\ar[r]^-{H_F}&S
}\]
is homotopy equivalent to $H$.

\begin{lemma}
The mapping class group $\cM(S_{g,n};\p)$ acts on $E\Fatbgn$ with quotient $\Fatbgn$.
\end{lemma}
\begin{proof}
$\cM=\cM(S_{g,n};\p)$ acts on $E\Fatbgn$ by composition with the marking.
It suffices to show that $\cM$ acts transitively on the markings of a
fixed bordered fat graph $\Gb$. Hence given two such markings 
$[H_i]$, we will construct $[\phi]\in \cM(S;\p)$ such that $[\phi\cdot H_1]=[H_2]$.

We first construct an homeomorphism $f \in \Homeo(S;\p)$ preserving the boundary pointwise such that $f\cdot H_1=H_2$. On
the image of $H_1:\Gb\hookrightarrow S$, define $f=H_2\cdot H_1^{-1}$.If $S_j$ denotes $S$ cut along the image of $H_j$ then each connected component $T_{ji}$ of $S_j$ is a polygon whose boundary contains exactly
one boundary component $\p_i S$ of $S$. The remainder of $\p T_{ji}$ is identified with the image of the \ith boundary cycle $\om_i$ inside
$H_j(\Gb)$. We have an orientation-preserving homeomorphism on the boundary of the
polygons, which extends to an orientation-preserving homeomorphism on the interior of each $T_{1i}$.

We approximate $f$ by a diffeomorphism. By a theorem of
Nielsen \cite{Nielsen},   $\Diff(S;\p)\hookrightarrow \Homeo(S;\p)$ is a homotopy equivalence and 
\[ \cM(S;\p)= \pi_0\Diff(S;\p)\cong \pi_0\Homeo(S;\p).\]
In particular, there is a diffeomorphism $\phi$ isotopic to $f$ and
\[[\phi\cdot H_1]=[f\cdot H_1]=[H_2].\]
\end{proof}

We are left to show that each category $E\Fatbgn$ realizes to a contractible space. For this we will
show that it is isomorphic to the arc complex constructed by Harer in \cite{Harer}. Let $y_i$ be new marked points on each of the boundary components of $S$ and let  
\[\Delta = \{y_1\ldots y_n\}.\]
An \textbf{$\Delta$-arc} is an embedded path or loop in $S$ starting and ending in $\Delta$ and otherwise disjoint from $\p S$.
An \textbf{arc-system} $[\alpha_0\ldots\alpha_k]$ is an isotopy class of $\Delta$-arcs such that
\begin{enumerate}
\item The intersections $\alpha_i\cap \alpha_j$ is contained in $\Delta$
\item For each component $G$ of the surface obtained by splitting $S$ along the $\alpha_i$'s, the Euler characteristic of the double of $\p
  G\setminus \Delta$ is negative.
\end{enumerate}
Note that by the second condition, no arc of an arc-simplex is null-homotopy relative to $\Delta$, no two arcs are
homotopic relative to $\Delta$ and no arc is homotopic to $\p_i S$ relative to $\Delta$.  
An arc-system $[\alpha_0\ldots \alpha_k]$ \textbf{fills} $F$ if all connected components of the $S\setminus\{\alpha_i\}$
are disks. Let $A_0(S)$ be the poset of filling arc-systems in $S$.

\begin{theorem}[\cite{Harer}]
There is an $M(S;\p)$-equivariant homotopy equivalence
\[|A_0(S_{g,n})| \map \cT_{g,n}\]
where $\cT_{g,n}$ denotes the Teichm\"uller space of $S_{g,n}$.
\end{theorem}

\begin{figure}
\begin{center}
\mbox{
\subfigure{\epsfig{file=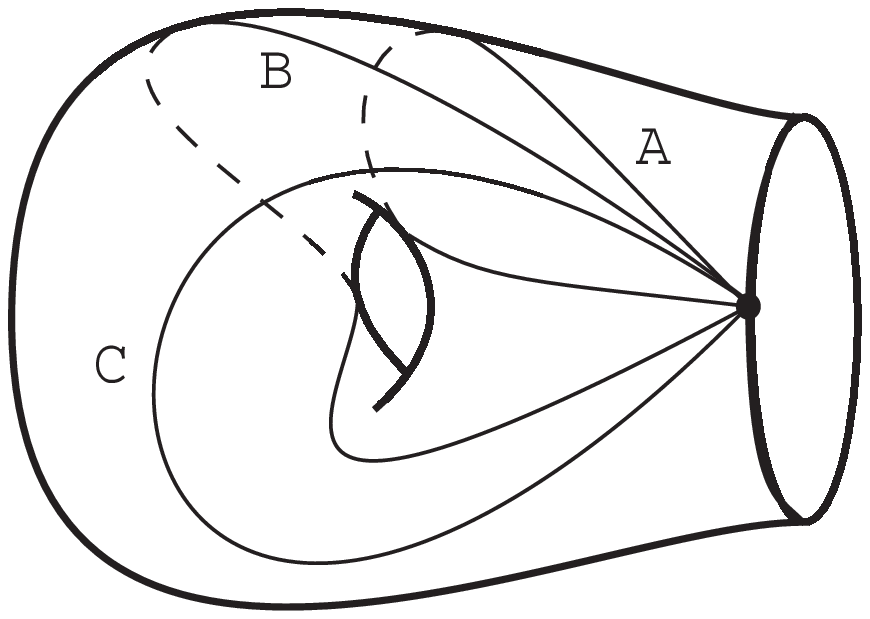, height=120pt}}\qquad
\subfigure{\epsfig{file=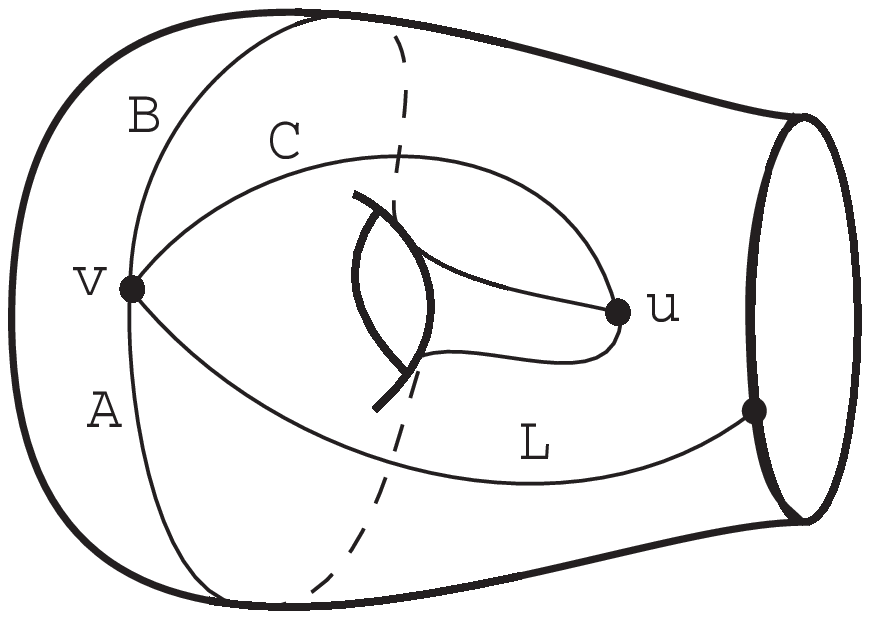, height=120pt}}}
\caption{An arc-system and its dual marked bordered fat graph}\label{fig:dual}
\end{center}
\end{figure}

\begin{lemma}
There is an $\cM(S;\p)$-equivariant isomorphisms
\[\Psi: A_0(S_{g,n}) \maplu{\cong} E\Fatbgn.\]
\end{lemma}
\begin{proof}
Given a filling arc-system $[\alpha_0,\ldots, \alpha_k]$, we build its dual
graph $G$ as follows. $G$ has one vertex $v_T$ for 
each component $T$ of  
\[S_0= S\setminus\{\alpha_0,\ldots,\alpha_k\}\]
and a vertex $v_i$ for each boundary component $\p_i S$. $G$ has an edge $E_j$
between $v_{T_1}$ and $v_{T_2}$ for each arc $\alpha_j$ touching both $T_1$ and
$T_2$ and an edge $\{L_i,\ol{L}_i\}$ between $v_i$ and the vertex corresponding
to the component $T_i$ of $S_0$ whose boundary contains $\p_i S$. 

Since the arc-system is filling, each component $T$ of $S_0$ is a polygon and
the orientation on $S$ gives a cyclic ordering $\sigma_v$ of the edges incident to the
vertices $v_T$ of $G$. The \ith boundary cycle of $\Gb=(G,\sigma)$ corresponds
to the arcs leaving $y_i$ and contains only one leaf $v_i$.
The bordered fat graph $\Gb$ has a natural marking $[H]$
where a vertex $v_T$ is sent in the region $T$ and where an 
edge $E_j$ intersects only the arc $\alpha_j$. 

Two isotopic sets of $\Delta$-arcs divide $S$ is an isomorphic 
fashion. Hence the isomorphism class of the marked bordered fat $[\G,[H]]$
depends only on the arc-system $[\alpha_0,\ldots\alpha_k]$.  
A similar construction also gives an arc-system $\gamma$ from a marked fat graph
$[\Gb,[H]]$. 

$\Psi$ is an isomorphism of poset. Assume
\[A= [\alpha_1,\ldots,\alpha_n]\subset B=[\alpha_1,\ldots, \alpha_n,
  \beta_1,\ldots \beta_m].\] 
The surface $S\setminus (\cup \alpha_i)$ is obtained from $S\setminus 
((\cup\alpha_i)\cup(\cup\beta_i))$ by gluing, for each $\beta_i$,  the two
components $T_i$ and $\tilT_{i}$ touching 
$\beta_i$ along $\beta_i$. In the dual graph, this corresponds to collapsing the
edge $E_{\beta_i}$. 
Similarly $\Psi(A)\leq \Psi(B)$ imply that $A\leq B$.
\end{proof}

\begin{proof}[Proof of theorem \ref{thm:M(S;p)}.] 
The action of the bordered mapping class group on $|E\Fatbgn|$ is free with quotient $|\Fatbgn|$.
If $(g,n)\neq (0,2)$ then by the previous lemma
\[|E\Fatbgn|\simeq \cT_{g,n}\simeq *.\]
\begin{figure}
\begin{center}
\mbox{
\subfigure[$G_0^b$]{\epsfig{file=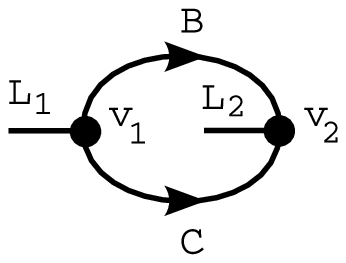, height=60pt}}\qquad
\subfigure[$G_1^b$]{\epsfig{file=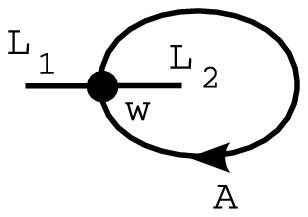, height=50pt}}}
\caption{The objects of $\Fat^b_{0,2}$}\label{fig:finb12}
\end{center}
\end{figure}
If $(g,n)=(0,2)$ then $S$ is a cylinder and $\cM(S;\p)=\bZ$ is generated by a Dehn twist along a meridian. The category
$\Fat^b_{0,2}$ has two objects $[G_0^b]$ and $[G_1^b]$ which are shown in figure
\ref{fig:finb12}. There are two non-identity morphisms 
$$\xymatrix{
[G_1^b]\ar@/^/[r]^B\ar@/_/[r]_C&[G_0^b]
}$$
and therefore
\[|\Fat^b_{0,2}|\simeq S^1\simeq B\cM(S_{0,2};\p).\]
\end{proof}

\begin{remark}\label{rmk:metric}
Using the argument of Igusa \cite[section 8.1]{Igu00}, 
one can show that the geometric realization of $\Fatb$ corresponds to a
space of isometric classes of metric bordered fat graphs. A vertex $[\Gb]$ represents a
metric graph whose edges have the same length. A morphism $[\Gb]\smap[\Gb/F]$ brings the length edges of $F$ to zero in a linear way.
\end{remark}

\section{A CW structure and its cellular graph complex.}
\label{sec:CW} 

Following ideas of Penner \cite{Penner_complex}, we
construct a chain complex $\GC^*$ generated by equivalence classes of bordered fat graphs. Note that $\GC^*$ computes
the \emph{integral} cohomology of the bordered mapping class group, where Penner's complex gives only the rational cohomology of the punctured ones. 

\subsection{CW structure on $|\Fatb|$.}
\label{sub:CW}
Using the combinatorics of the bordered fat graphs, we construct a filtration on $\Fatb$ and show that it induces a CW
structure on the geometric realization $|\Fatb|$ with exactly one cell per isomorphism class of bordered fat graphs.
\begin{definition}
The {\bf dimension} of a bordered fat graph $\Gb$ is equal to
\begin{equation*}
\codim(\Gb) =\sum_{v}\big(val(v)-3\big) 
\end{equation*}
where the sum is taken over all interior vertices $v$ and where $val(v)$ is the number of half edges incident to $v$. 
\end{definition}

Recall that the interior vertices have valence at least three. A bordered fat graph has dimension zero if all its interior (non-leaf)
vertices have valence exactly three, such a fat graph will be called \emph{essentially trivalent}. A morphism of bordered fat graphs increases the dimension by the
number of edges it collapses. 

Let $\cF^p\Fatb$ denote the full subcategory of $\Fatb$ containing all bordered fat graphs of dimension at
most $p$.

\begin{theorem}\label{thm:fp/fp-1}
The geometric realization of the filtration
\[ \cF^0\Fatb\subset\cdots\subset \cF^p\Fatb\subset\cdots \subset \Fatb\] 
is the skeleton of a CW-structure on $|\Fatb|$ which has a $p$-cell for each isomorphism classes of bordered fat graph of dimension $p$. 
\end{theorem}

\begin{remark}\label{rmk:cell}
The $p$-cell of $[\G^b_p]$ contains exactly the p-simplices
\begin{equation*}
[\G^b_0]\smap\ldots\smap[\G^b_{p-1}]\smap [\G^b_{p}]
\end{equation*}
where each $\G^b_i$ has dimension $i$.
\end{remark}

A morphism $[\varphi:\tilGb\smap\Gb]$ corresponds to a choice of a tree $T_v$ to replace each vertex of $\Gb$.
Hence the category $\Fatb/[\G]$ is a product of categories of planar trees as follows.

We define $\cT^k$ to be the category whose objects are isomorphism classes of
combinatorial planar trees $[T]$ with $k$ fixed cyclically-ordered leaves and
whose interior vertices are at least trivalent. For each forest $F$ of $T$
containing only interior edges, $\cT^k$ has a morphism  
\[\psi_{[F]}:[T]\map [T/F]\]
where $T/F$ is the planar tree obtained from $T$ by collapsing each the edges of
$F$. In particular the morphisms preserve the $k$ leaves and their cyclic
ordering. Let $\cT^k_1$ denote the full subcategory of $\cT^k$ containing all
trees with at least one interior edge.  

\begin{lemma}[Stasheff \cite{Stasheff}]\label{lem:assoc}
The geometric realization of the category $\cT^k$ is homeomorphic to a cube
$I^{k-3}$ and the subcategory $\cT_0^k$ realizes to its boundary. 
\end{lemma}

\begin{proof}[Sketch of Proof.]
The category $\cT^k$ realizes to a space $X_k$ of isometric classes of metric
plannar trees with cyclically ordered leaves $\{1,2,\dots k\}$. In $X_k$, only
interior edges are given lengths. 
The claim is proved by induction on $k$. $X_3$ is a single point. From each
object of $X_{k-1}$, one builds an object of $X_{k}$ by attaching the \xth{k} leaf
between the last and the first. Hence $X_{k+1}$ is an interval bundle over
$X_{k}$.
\end{proof}

\begin{remark}
In \cite{Stasheff} Stasheff constructed the space $K_{k}\cong X_{k+1}$ without
the use of trees. He defined the $K_k$ inductively using all meaningful way of
inserting a pair of parentheses in the word $x_1\ldots x_{k}$.  
\end{remark}

\begin{lemma}\label{lem:fattrees}
For any isomorphism class of bordered fat graphs $[\Gb]$ the category $\Fatb/[\Gb]$ is isomorphic to the product category
\[\Fatb/[\Gb]\esp{\cong}\prod_{v} \cT^{val(v)}\]
where the product is taken over all interior vertices. 
\end{lemma}

\begin{proof}
An object in $\Fatb/[\Gb]$ is a morphism $\psi_{[\tilGb,F]}$ of $\Fatb$ such that $\tilGb/F$ is
isomorphic to $\Gb$. A morphism of $\Fatb/[\Gb]$ between $\psi_{[\G^b_1,F_1]}$ and $\psi_{[\G^b_2,F_2]}$ is a morphism
$\psi_{[\G^b_1,E]}$ of $\Fatb$ so that the following diagram commute.
\[\xymatrix{
[\G^b_1]\ar[rd]_{\psi_{[\G^b_1,F_1]}}\ar[rr]^{\psi_{[\G^b_1,E]}}&&[\G^b_2]\ar[ld]^{\psi_{[\G^b_2,F_2]}}\\
&[\Gb].
}\]

For an object $\psi_{[\tilGb,F]}$, cut the edges of $\tilGb$ that are not in $F$ in two. After this procedures $\tilGb$ becomes a disjoint union of trees, one $T^v$ for 
each vertex $v$of $\Gb$. The tree $T^v$ has one leaf for each half-edge adjacent to $v$ in $\Gb$. A morphism $\psi_{[\G^b_1,E]}$ collapses a subforest $E$ of $F_1$. At each vertex $v$, it collapses the subforests $E\cap T^v_1$ of the tree $T^v_1$ collapsed by $\psi_{[\G^b_1,F_1]}$. Using this we get an isomorphism
\begin{eqnarray*}
\Fatb/[\Gb]&\map&\prod_{v} \cT^{val(v)}\\
\psi_{[\G^b,F]}&\maps&\big([T^v]\big)_{v}\\
\big([\G^b_1,F_1]\smaplu{[E]}[\G^b_2,F_2]\big)&\maps&\big([E\cap\tilT_v^1]: [T^v_1] \map [T^v_2]\big).
\end{eqnarray*}
\end{proof}

\begin{proof}[Proof of theorem \ref{thm:fp/fp-1}]
The boundary of the $p$-cell $|\Fatb/[\Gb]|$ is the geometric realization of the full subcategory of $\Fatb/[\G]$ containing the
objects $\psi_{[\tilGb, F]}$ where $F$ is not empty. In particular, such a $\tilGb$ has dimension
strictly less than $p$.
Let 
\[X = \Big(|\cF^{p-1}\Fatb| \cup \bigcup_{\codim\Gb=p}|\Fatb/\Gb|\Big)\]
be the space obtained by attaching the boundary of each $|\Fatb/\Gb|$ along the forgetful functor
\[\p(\Fatb/\Gb) \map \cF^{p-1}\Fatb.\]

We have a natural map $\pi:X\map  |\cF^p\Fatb|$. Fix a simplex
\[\alpha =\big([\G^b_0] \map\cdots\map[\G^b_k]\big).\]
of $\cF^p\Fatb$. If $\codim(\G^b_k)<p$ then $\alpha$ is a simplex of $\cF^{p-1}\Fatb$ and appears once in $X$. Otherwise $\codim(\G^b_k)=p$
and is in the image of
\[[\G^b_0] \map\cdots\map[\G^b_k]\big)\maplu{Id} [\G^b_k] \in \Fatb/[\G^b_k].\]
Since all non-identity morphisms in $\Fatb$ collapse at least an edge and increase dimension, $\alpha$ is the image
of a single simplex of $X$. Hence $\pi$ is an homeomorphism.
\end{proof}

\subsection{The bordered graph complex.}
\label{sub:complex}

In this section, we describe the bordered graph complex which is the cellular cochain complex associated to the CW-structure on
$|\Fatb|$ introduced in section the previous section. Most of this section is a translation to our setting of idea found in \cite{Penner_complex} and \cite{ConVog}.

\begin{definition}
An {\bf orientation} on a graph $G$ is a unit vector in 
\[\cO(G)= det(\bR V_G) \tens det(\bR H_G)\]
where $V_G$ is the set of vertices of $G$ and where $H_G$ is its set of half-edges.
\end{definition}

\begin{remark}\label{rem:orient}
An orientation $o$ on $\Gb$ corresponds to an equivalence class of an ordering of the vertices and an orientation of
each edge. This class is obtained by reordering the terms in $o$ as follows.
\[o= (v_1\wedge \cdots\wedge v_n) \tens (A_1\wedge\ol{A}_1\wedge A_2\wedge \ol{A}_2\cdots)\]
A bordered fat graph $\G$ of dimension zero has a natural orientation $o_{\nat}(\G)$ induced by its fat structure.
Choosing an ordering $v_1\ldots v_k$ of the (non-leaf) vertices of $\G$. Choose a total ordering $a_1^v\cdots a^v_k$ compatible with the
cyclic ordering $\sigma(v)$. Define
\begin{equation}\label{eq:natorient}
o_{\nat} = (v_1\wedge \cdots\wedge v_k)\tens (a_1^{v_1}\wedge\cdots\wedge
a^{v_1}_{k_1} \wedge a_1^{v_2}\wedge \cdots)
\end{equation}
Since the valence of each vertex is odd, this orientation is independent of the choices made.
\end{remark}

\begin{proposition}\label{prop:orient}
Given a bordered fat graph $\Gb$ and an edge $e$ which is not a cycle, there is an isomorphism 
\[ \Phi_e : \cO(\Gb)\maplu{\cong}\cO(\Gb/e)\]
between the orientations of $\Gb$ and the ones of the bordered fat graph $\Gb/e$  obtained from $\Gb$ by collapsing $e$
\end{proposition}
\begin{proof}
Let $e=\{A,\ol{A}\}$ and let $o$ be an orientation of $\Gb$. Rewrite $o$ so that
\[ o = (s(A)\wedge s(\ol{A}) \wedge u_1\wedge \ldots)\tens (A\wedge \ol{A} \wedge h_1\wedge \ldots).\]
Denote by $\widetilde{u}$ the image of $e$ in $\Gb/e$ and define 
\[ \Phi_e(o) = o_e = (\widetilde{u}\wedge u_1\wedge \ldots)\tens (h_1\wedge \ldots).\]
\end{proof}

Let $\GC^p$ be the free abelian group generated by the equivalence classes $[\Gb,o]$ of oriented bordered fat graphs
modulo the relation 
\[-[\Gb,o] = [\Gb,-o].\]
We define a coboundary map $\delta:\GC^p\smap \GC^{p+1}$  
\[\delta \big([\G^b_p,o]\big) = \sum [\G^b_{p}/e, \Phi_e(o)]\]
where we sum over all edges $e$ of $\G^b_p$ which are not cycles. 
Similarly we define $\GC_p =\GC^p$ and a boundary map $\p:\GC_p\smap \GC_{p-1}$
\[\p [\G^b_p,o] = \sum [\G^b_{p-1},\Phi_e^{-1}(o)]\]
where we sum over isomorphism classes $[\G^b_{p-1},e]$ such that $ \G^b_p$ is obtained from $\G^b_{p-1}$ by collapsing $e$.

\begin{theorem}\label{thm:K chain complexes}
\begin{eqnarray*}
H^*(\GC^*;\delta)&\cong& H^*(|\Fatb|;\bZ)\esp\cong \bigoplus\ H^*(\cM(S;\p);\bZ)\\
H_*(\GC_*;\p) &\cong& H_*(|\Fatb|;\bZ) \esp\cong \bigoplus\ H_*(\cM(S;\p);\bZ)
\end{eqnarray*}
\end{theorem}

\begin{remark}
Bordered fat graphs have no automorphism. Because of this, the bordered graph complex computes the \emph{integral} (co)homology of the mapping class groups, eventhough the original graph complex computed only the \emph{rational} one.
\end{remark}

To prove theorem \ref{thm:K chain complexes}, we will show that an orientation on $\Gb$ corresponds to a choice of
compatible orientations for each $p$-simplex of the $p$-cell at $[\Gb]$. Given an oriented bordered fat graph $(\Gb,o)$ of dimension $p$ and given a $p$-simplex
\[\alpha := \big(\G^b_0\maplu{\{e_1\}}\G^b_0/e_1=\G^b_1\maplu{\{e_2\}} \ldots\maplu{\{e_p\}} \Gb\big),\]
we define $\epsilon(\alpha,o_{\G_p})$ to be $\pm 1$ depending on whether the
orientation on $\G^b_p$ induced by $\alpha$ and $o_{\nat}(\G^b_0)$ agrees or disagrees with $o$. Hence
\[\Phi_{e_p}\cdots \Phi_{e_1} (o_{\nat}(\G^b_0)) =  \epsilon(\alpha,o)\, o.\]

\begin{lemma} \label{lem:orient pcells}
Given an oriented bordered fat graph $(\G^b_p,o)$ of dimension $p$, hte chain
\begin{equation}\label{eqn:sphere}
\beta = \sum_{\alpha = [\G_0\smap\cdots\smap\G_p] } \epsilon(\alpha, o)\ [\G^b_0\smaplu{\{e_1\}}\cdots\smaplu{\{e_p\}}\G^b_p]\maplu{Id} [\G^b_p]
\end{equation}
is an orientation of $|\Fatb/[\Gb]|$ relative to its boundary.
\end{lemma} 

\begin{proof}
Since each $p$-simplices of $\Fatb/[\Gb]$ appears exactly once in $\beta$, it suffices to show that the boundary of \eqref{eqn:sphere} lands in the boundary of $\Fatb/[\Gb]$.
Modulo simplices of the boundary of $\Fatb/[\Gb]$, we have
\begin{equation*}
\begin{split}
\p(\beta) {\equiv} \sum_{\alpha}\ o(\alpha)\ \Big(\sum_{i=1}^{p-1}\ (-1)^i\,
  &[\G^b_0\smap\cdots \xymatrix{\G^b_{i-1}\ar[r]^{\{e_i,e_{i+1}\}}&\ \G^b_{i+1}}\cdots\smap\G^b_p]\maplu{Id}[\Gb]\\
  &+ [\G^b_1\smaplu{\{e_1\}}\cdots\smaplu{\{e_p\}}\G^b_p]\maplu{Id}[\Gb]\Big)
\end{split}
\end{equation*}
A simplex $[\cdots\, \G^b_{i-1}\smap\G^b_{i+1}\cdots]$ appears in the boundary of two simplices of $\beta$ in
$\p(\beta)$, namely
\[\alpha_1 = [\cdots \G^b_{i-1}\maplu{\{e_{i}\}}\G^b_{i}\maplu{\{e_{i+1}\}}\G^b_{i+1}\cdots] \qquad \alpha_2=[\cdots
  \G^b_{i-1}\maplu{\{e_{i+1}\}}\tilG^b_{i}\maplu{\{e_{i}\}}\G^b_{i+1}\cdots].\]
It suffices to show that the induced orientations are opposite.
Let $\tilde o$ be the orientation on $\G^b_{i-1}$ induced by the natural one on $\G^b_0$. Assume $e_i=\{A,\ol{A}\}$ and
$e_{i+1}=\{B,\ol{B}\}$ do not share a vertex. If 
\[ \tilde o = (s(A)\wedge s(\ol{A})\wedge s(B) \wedge s(\ol{B})\wedge v_1 \cdots)\tens(A\wedge \ol{A}\wedge B\wedge\ol{B}\wedge
h_1\cdots)\]
then one gets
\begin{eqnarray*}
\Phi_{e_{i+1}}\Phi_{e_i} (\tilde o) &=& (\beta\wedge \alpha\wedge v_1\cdots)\tens(h_1\cdots)\\
\Phi_{e_i}\Phi_{e_{i+1}} (\tilde o) &=& (\alpha\wedge\beta\wedge v_1\cdots)\tens(h_1\cdots)
\end{eqnarray*}
where $\alpha$ is the image of $e_i$ and $\beta$ the one of $e_{i+1}$.
Hence, 
\[\Phi_{e_p}\cdots\Phi_{e_{i+1}}\Phi_{e_{i}}\cdots\Phi_{e_1}(o_{\nat}) = -
\Phi_{e_p}\cdots\Phi_{e_{i}}\Phi_{e_{i+1}}\cdots\Phi_{e_1}(o_{\nat}).\]
The case where $e_i$ and $e_{i+1}$ share a vertex is similar and will be omitted.

The simplex $[\G^b_1\smap\cdots]$ also appears in the boundary of two simplices in $\p(\beta)$. A
similar reasonning gives that the two induced orientations will cancel.
\end{proof}

\begin{remark} \label{rmk:chain}
In the previous proof, we have shown that
\[ \Phi_e \Phi_f = - \Phi_f\Phi_e.\]
This implies that $\p^2\equiv0$ and that $\GC$ is a chain complex.
\end{remark}

\begin{proof}[Proof of theorem \ref{thm:K chain complexes}.]
By section \ref{sub:CW}, we have a CW structure on $|\Fatb|$ whose cells come from $\Fatb/[\Gb]$. The previous lemma
therefore identifies these cells in $C^{\mathit{simp}}_p$. To show that $\GC_*$ is the cell complex, it suffices to
prove that the following map 
\begin{alignat*}{10}
\Phi_p\esp{:}&& \GC_p &\map C^{\mathit{simp}}_p\\
&&[\G^b_p,o]&\maps  (-1)^{(p+1)p/2} \sum_{\alpha = [\G_0\smap\cdots\smap\G_p] } o(\alpha, o)\ [\G^b_0\smap\cdots\smap\G^b_p].
\end{alignat*}
is a chain map.

By the argument lemma \ref{lem:orient pcells}, we know that
\[\begin{split}
\p \sum_{\alpha = [\G_0\smap\cdots\smap\G_p] } o(\alpha, o) \alpha 
&\esp{=} \sum_{\alpha = [\G_0\smap\cdots\smap\G_p] } o(\alpha, o)\, (-1)^p\,  [\G_0\smap\cdots\smap\G_{p-1}]\\
&\esp{=} (-1)^p\sum_{[(\G_{p-1},\tilde o)\smap(\G_p,o)]} \sum_{\beta=[\G_0\smap\cdots\smap\G_{p-1}]} o(\beta, \tilde o)\beta
\end{split}\]
Each term in the last sum  is a representative for a $p-1$ cell of $\Fatb$ and the signs work out so that
\[ \p\cdot\Phi_p[\G_p] = \Phi_{p-1} \delta [\G_p].\]
This proves the second statement. The first one is simply its dual.
\end{proof}

\subsection{Computations.}

The bordered  graph complex can be separated into smaller cochain complexes each
computing the cohomology of a single connected component of $|\Fatb|$ and
therefore of a single mapping class group. In this section we compute
explicitely the homology and the cohomology of $\cM(S_{1,1};\p)$. We also 
give the homology of $\cM(S_{1,2};\p)$ and $\cM(S_{2,1};\p)$ which were obtained
with the help of a computer.

\label{sub:homS11}

\begin{figure}
\begin{center}
\mbox{
\subfigure[$\G^b_2$]{\epsfig{file=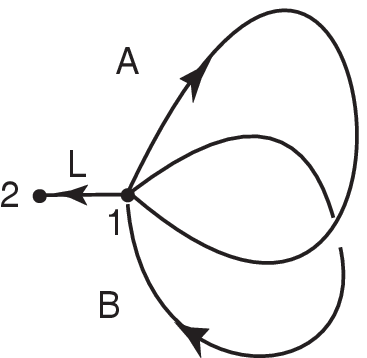,height=65pt}}\quad
\subfigure[$\G^b_1$]{\epsfig{file=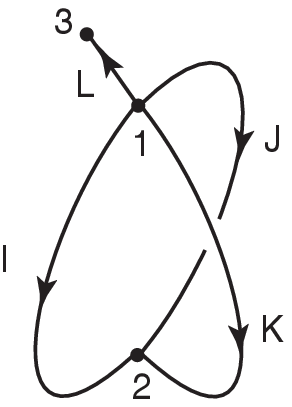,height=75pt}}\qquad
\subfigure[$\tilG^b_1$]{\epsfig{file=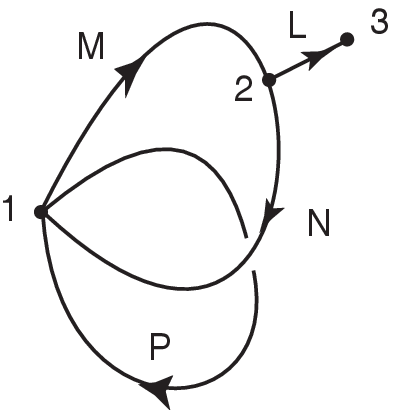,height=70pt}}\quad
\subfigure[$\G^b_0$]{\epsfig{file=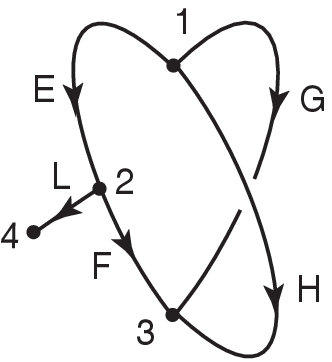,height=70pt}}}
\end{center}\caption{The objects of $\Fat^b_{1,1}$}\label{fig:objects}\label{fig:genS11}\label{fig:s11}
\end{figure}

\begin{proposition}\label{prop:S11} The homology and cohomology of the mapping class group  of
  the onced-bordered torus $S_{1,1}$ is as follows.
\[H^*(\cM(S_{1,1};\p))\cong H_*(\cM(S_{1,1};\p),\bZ)= \begin{cases}\bZ & *=0,1\\0 & *\geq 2.\end{cases}\]
The homology is generated by the cells associated to the bordered fat graphs $\G_0^b$
and $\tilG^b_1$ of figure \ref{fig:objects}. 
\end{proposition}

\begin{proof}
There are four isomorphism classes of bordered fat graph in $\Fat^b_{1,1}$, all of which are shown
in figure \ref{fig:objects}. 
Using remark \ref{rem:orient}, we give each of these bordered fat graphs the orientation
induced by the ordering of the vertices and the orientation of the edges
depicted in figure \ref{fig:objects}. For example, 
\[o_{\G^b_2} = (1\wedge 2)\tens (A\wedge \ol{A}\wedge B\wedge \ol{B}\wedge L\wedge \ol{L}).\]
The coboundary of each of these oriented bordered fat graphs is
\[\delta([\G_0^b]) =0\qquad \delta([\G_1^b]) = [\G^b_2]\qquad
\delta([\tilG^b_1]) =2*[\G^b_2]\]
which gives the cohomological result. 
We will calculate the coboundary of $\tilG^b_1$ explicitly. By
definition,
\[\delta\big([\tilG^b_1,o]\big) = [\tilG^b_1/M, \Phi_M(o)] + [\tilG^b_1/N,\Phi_N(o)]\]
Both $\tilG^b_1/M$ and $\tilG^b_1/N$ are isomorphic to
$\G^b_2$. Collapsing $M$ gives
\[\begin{aligned}
o_{\tilG^b_1} 
&= (1\wedge 2\wedge 3)\tens (M\wedge\ol{M}\wedge N\wedge\ol{N}\wedge
P\wedge\ol{P} \wedge L\wedge\ol{L})\\
&\mapsto (\{1,2\}\wedge 3 )\tens (N\wedge \ol{N} \wedge P\wedge\ol{P} \wedge L\wedge\ol{L})\\
&\cong  (1\wedge 2)\tens (A \wedge \ol{A} \wedge B\wedge \ol{B}) = o_{\G^b_2}
\end{aligned} \]
Similarly, when collapsing $N$,
\[\begin{aligned}
o_{\tilG^b_1}  
&= - (2\wedge 1\wedge 3)\tens (M\wedge\ol{M}\wedge N\wedge\ol{N}\wedge
P\wedge\ol{P} \wedge L\wedge\ol{L})\\
&\mapsto - (\{2,1\} \wedge 3)\tens (M\wedge\ol{M}\wedge P\wedge\ol{P}
\wedge L\wedge\ol{L})\\
&\cong - (1\wedge 2)\tens (B\wedge \ol{B}\wedge \ol{A}\wedge A \wedge
L\wedge \ol{L}) = o_{\G^b_2}
\end{aligned}\] 
Hence, as claimed, $\delta([\tilG^b_1])= 2[\G^b_2]$.
\end{proof}

\begin{remark}\label{rem:s11}
The mapping class group of the surface $S_{1,1}$ has the following presentation  
\[ \cM(S_{1,1};\p)= \big< \alpha, \beta \btab \alpha \beta\alpha =
\beta \alpha\beta\big>.\]
where $\alpha$ and $\beta$ are Dehn twists around the meridian and the longitude
respectively. $\cM(S_{1,1};\p)$ is therefore isomorphic to the braid group on three
strands. By a theorem of Quillen \cite[p.84]{Milnor} the complement of
the trefoil knot is a classifying space for this braid group. It follows that $B\cM(S_{1,1};\p)$ is an homology sphere. 
\end{remark}

Although the surfaces $S_{1,2}$ and $S_{2,1}$ are still relatively
simple, the categories $\Fat^b_{1,2}$ and  $\Fat^b_{2,1}$ have over a
thousand objects. Their associated bordered graph complexes are
therefore too big to study manually. 

To construct the appropriate bordered graph complex, we use a computer
program. We first find all isomorphism classes of bordered fat graphs of the
desired type and then construct the bordered graph complex by computing the
boundary maps of the graph complex. We then input this bordered graph complex
into the computer algebra magma and get the following results. 

\begin{proposition} \label{prop:s21}\label{prop:s12}
The homology groups of the bordered mapping class group of $S_{1,2}$ are
\[H_*(\cM(S_{1,2};\p))\cong\begin{cases}
\bZ &*=0\\
\bZ\oplus\bZ & *=1\\
\bZ/2\oplus \bZ & *=2\\
\bZ/2 & *=3\\
0 & *\geq 4.
\end{cases}\]
The homology groups of the bordered mapping class group of $S_{2,1}$ are 
\[H_*(\cM(S_{2,1};\p)) \cong\begin{cases}
\bZ & *=0\\
\bZ/10& *=1\\
\bZ/2 & *=2\\
\bZ\oplus \bZ/2 & *=3\\
\bZ/6 &*=4\\
0 &*\geq 5.
\end{cases}\]
\end{proposition}

\begin{figure}
\begin{center}
\mbox{
\subfigure{\epsfig{file=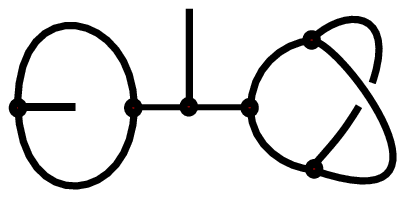, width=60pt}}\quad
\subfigure{\epsfig{file=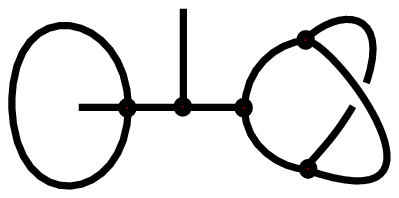, width=60pt}}\quad
\subfigure{\epsfig{file=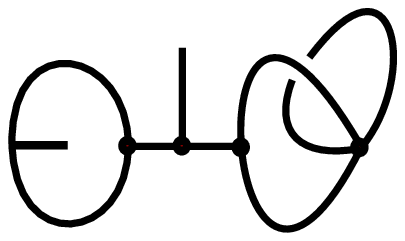, width=60pt}}\quad
\subfigure{\epsfig{file=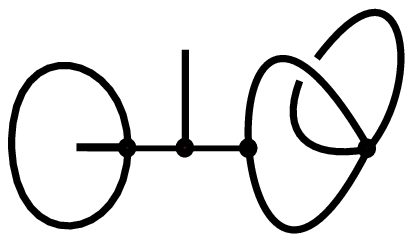, width=60pt}}\quad
\subfigure{\epsfig{file=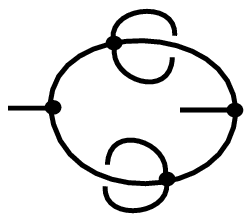,
    width=45pt}}}
\caption{Generators for $H_k(\cM(S_{1,2};\p))$ for $k=0,1,2$}\label{fig:gen HS12}
\end{center}
\end{figure}

\begin{figure}
\begin{center}
\mbox{
\subfigure{\epsfig{file=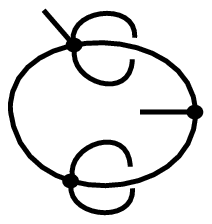, width=50pt}}\qquad
\subfigure{\epsfig{file=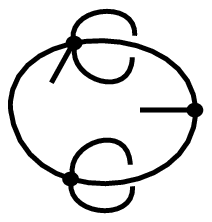, width=50pt}}\qquad
\subfigure{\epsfig{file=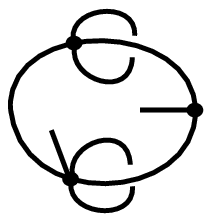, width=50pt}}\qquad
\subfigure{\epsfig{file=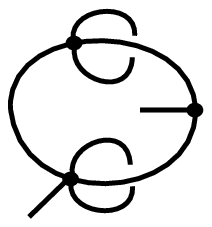, width=50pt}}}
\caption{Generator for $H_4(\cM(S_{1,2};\p))$}\label{fig:gen H4S12}
\end{center}
\end{figure}

\begin{figure}
\begin{center}
\mbox{
\subfigure{\epsfig{file=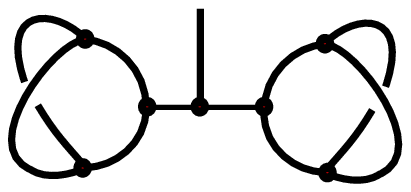, width=90pt}}\qquad
\subfigure{\epsfig{file=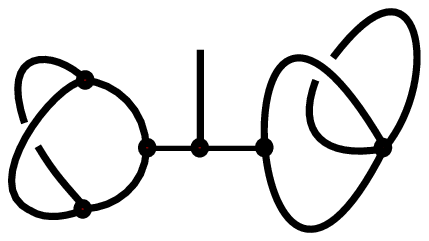, width=90pt}}\qquad
\subfigure{\epsfig{file=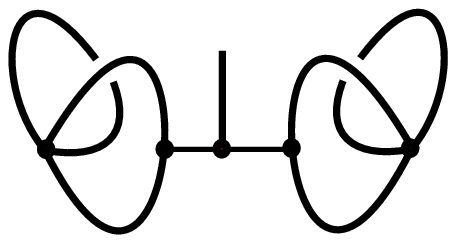, width=90pt}}}
\caption{Generators for $H_k(\cM(S_{2,1};\p))$ for $k=0,1,2$ }\label{fig:genS21} \label{fig:s21}
\end{center}
\end{figure}

The computer algebra magma give generators for these homology groups. However, they involve many cells and they are unenlightning. The generators of figures \ref{fig:s21}, \ref{fig:gen HS12} and \ref{fig:gen H4S12} are obtained by using the circle action of section \ref{sec:DX} and Miller's operad of section \ref{sec:Miller}.

\begin{remark}
Some of these results were already known. Using a spectral sequence argument Ehrenfried in \cite{ehr}  got the same result for the homology of $\cM(S_{2,1};\p))$.
\end{remark}

We have also constructed the restriction of the bordered graph complex to
subcategories $\Fat^b_{3,1}$, $\Fat^b_{2,2}$ and $\Fat^b_{1,3}$. However these bordered 
graph complexes have over a million objects.

\section{The integral homology of the unbordered mapping class group.}
Let the unbordered mapping class group of a surface $S$ 
\[\cM(S)=\pi_0\ \Diff(S;\p_1S;\dots;\p_nS)\]
be the group of isotopy classes of orientation-preserving self-diffeomorphism of $S$ which preserves
each boundary component as a set. The goal of this section is to find a combinatorial complex
computing the \emph{integral} homology of this punctured mapping class group. 

Consider a surface $S=S_{g,1}$ with a single boundary component. Its bordered mapping class group is
a $\bZ$-central extension of its unbordered one. In particular, there is a quasi-fibration
\[ B\cM(S;\p)\map B\cM(S)\map \CPi\]
which we now use to transfer our knowledge of the homology of the bordered mapping class group to
the homology of the unbordered one.

\subsection{An integral complex for the unbordered mapping class group.}
\label{sec:DX}
By twisting the bordered graph complex with a circle action on the bordered fat graphs, we get a double complex computing the integral homology of the unbordered mapping class group. 

For any $S^1$-space $X$, denote by $\Delta_X$ the map
\[\Delta_X:C_p(X)\map C_1(S^1)\tens C_p(X)\map C_{p+1}(S^1\times X)\map C_{p+1}(X)\]
obtained by multiplying a chain by the fundamental class of the circle. Let
\[D_X = C_*(X)\tens\bZ[u]\]
where $u$ has degree 2 and define a differentiation by
\begin{eqnarray*}
d_{tot}: C_*(X)\tens\bZ[u]&\map& C_*(X)\tens\bZ[u]\\
\alpha\tens u^{k} &\maps& \p\alpha \tens u^k + (-1)^{deg(\alpha)} \Delta_X(\alpha) \tens u^{k-1}.
\end{eqnarray*}

\begin{proposition}[{\cite[Lemma 5.1]{Jones}}] \label{prop:DX}
The chain complex $(D_{X}, d_{tot})$ computes the $S^1$-equivariant homology of $X$. It is natural in the $S^1$-space $X$. 
If $\tau\in H^2_{S^1}(X)$ is the Euler class then capping with $\tau$ lowers the degree of $u$ by $1$.
\end{proposition}
\begin{remark}
Jones uses cyclic sets to prove this statement. Here is a more geometric justification.
The chain $[\CP^1]$ of $\CPi$ lifts to 
\[\beta : D^2 \map \Si\]
which sends $z$ to $(z,a,0,\dots)$ for the appropriate $a\geq0$. For any chain $\alpha$ of $X$, we lift $\alpha\tens [\CP^1]$ to a chain 
\[[\alpha \tens\beta] \in C_*(X\underset{S^1}\times S^{\infty})\]
where $S^1$ acts diagonally. The boundary of this lifting gives
\[(d\alpha)\tens \beta +(-1)^{|\alpha|}\alpha\tens (d\beta)\]
which corresponds to
\[\big[(d\alpha)\tens\beta +(-1)^{|\alpha|} \alpha\tens S^1 \big] = \big[(d\alpha)\tens\beta +(-1)^{|\alpha|} \Delta_X(\alpha)\tens 1 \big]. \]
\end{remark}

To twist the bordered graph complex with the circle action, we construct a combinatorial
$S^1$-action on the bordered graph complex restricted to surfaces
with a single boundary component, which we shall denote by $\GC^1$.

Take a bordered fat graphs $\Gb$ with a unique boundary cycle
\[\om=(L\ol{L} A_1\ldots A_n).\]
Assume the leaf-edge $L$ is attached to a trivalent vertex $z$. We define $\Phi_i(\G^b)$ to be the
bordered fat graph obtained as follows. We remove the leaf-edge
$L$ and the leaf. We also remove $z$ which is now bivalent and join its two
remaining adjacent edges $A_1$ and $A_n$. Finally we attach a new
leaf-edge $M$ at the source of $A_i$ so that the 
new boundary cycle becomes 
\[(M\,\ol{M}\, A_i\ldots \{A_n, A_1\} \ldots A_{i-1}). \] 
For an orientation
\[ o_\G = (v_1\wedge \cdots \wedge v_{r}\wedge z\wedge l) \tens
(b_1\wedge\cdots \wedge A_n\wedge A_1\wedge \ol{A}_1 \wedge \ol{A}_n\wedge
L\wedge \ol{L})\]
on $\Gb$, we construct the following orientation on $\Phi_i(\Gb)$. 
\[\Phi_i(o_\G) = (v_1\wedge \cdots v_{r} \wedge l)\tens
(b_1\wedge\cdots\wedge \{A_n, A_1\} \wedge \ol{\{A_n,A_1\}} \wedge
M\wedge \ol{M})\]
Define a map
\begin{eqnarray*}
\zeta\ :\ \GC_p &\map& \GC_{p+1}\\
{}[\Gb,o_\G]&\maps& \begin{cases}
\sum [\Phi_i(\Gb), \Phi_i(o_\G)] & \val(z) = 3\\
0 &\val(z)\geq 4
\end{cases}
\end{eqnarray*}
where $z$ is again the vertex onto which the leaf-edge $L$ is attached.

\begin{lemma} \label{lem:Phi}\label{lem:circle}
The map $\zeta$ is a chain map on the bordered graph complex which
induces the map
\[H_p(\cM(S;\p))\map H_{p+1}(\cM(S;\p))\]
obtained by multiplying any cycle by the full circle.
\end{lemma}

We delay the proof of this lemma which is technical to give an immediate application of this lemma
which is a corollary to Jones' proposition.

\begin{corollary}\label{cor:ssCPi}
The chain complex
\[\GC^1_*\tens \bZ[u]\]
with derivative
\[ \p_{tot} \Big([\G^b_p]\tens u^k\Big) = \left(\p[\G^b_p]\right)\tens u^k + (-1)^p \big(\zeta([\G^b_p])\tens
u^{k-1}\big)\]
computes the \textbf{integral} homology of the unbordered mapping class groups $\cM(S_{g,1})$.
\end{corollary}

\begin{corollary}For any surface $S_{g,1}$, there is a spectral sequence with
\begin{eqnarray*}
E^2_{p,q} &=& \begin{cases}
  H_{p}(\GC_*) & q =2k\\
  0 &q=2k+1
\end{cases}\\
d_2 &=& [\zeta] : H_p(\GC) \map H_{p+1}(\GC)
\end{eqnarray*}
which converges to the homology of the punctured mapping class group.
\end{corollary}
\begin{remark}\label{rmk:ss}
Since we have a double complex, finding the derivatives of the spectral sequences is easy. Let $[\alpha]\in H_n(X)\cong E^2_{n,k}$ where $k$ is even. First, we find 
\[d_2 ([\alpha]) = \Delta_X(\alpha) \in E_{n+1,k-2}^2\]
If $[\alpha\tens u^k]$ survives to $E_{n,k}^4$, we find that
\[d_4([\alpha\tens u^k]) =\pm \Delta_X(\beta) \in E_{n+3,k-4}^4\]
where $\beta$ is any chain of $C_{n+2}(X)$ such that $\p\beta=\Delta_X(\alpha)$. And so on.
\end{remark}

The following is a result of B{\"o}digheimer and Tillmann  in \cite{Bod_Tillmann} and it follows directly from the preceeding corollary.
\begin{corollary}\label{cor:trivialS1}
The circle action 
\[H_p(\cM(S_{g,1};\p)) \map H_{p+1}(\cM(S_{g,1};\p))\]
on the homology of the bordered mapping class group of a surface is stably
trivial.
\end{corollary}
\begin{proof}
Consider the natural group homomorphisms
\[\cM(S_{g,1};\p)\map \cM(S_{g,1})\map \cM(S_{g,0})\]
where $S_{g,0}$ is a closed surface of genus $g$. The composition glues a disk
to the boundary component and by Harer stability, it induces an isomorphism in
homology in a range of dimension. The stable homology of the bordered mapping
class groups is therefore a summand of the stable homology of the unbordered
mapping class groups. In particular the stable part of the first column of our
spectral sequence survives and hence 
\[d_2=\Delta\ :\ E_{2,q}^2\map E_{0,q+1}^2\]
is zero in the stable range.
\end{proof}

\begin{proof}[Proof of lemma \ref{lem:Phi}]
Recall that by remark \ref{rmk:metric}, $|\Fatb|$ is homeomorphic to a space of metric bordered fat
graphs. Define a map
\[\widetilde\Phi\ :\ \bR\times |E\Fatb| \map |E\Fatb|\]
where $t\in \bR$ slides the leaf-edge $t$ times around the boundary
cycle of any marked metric bordered fat graph. This map is equivariant with regards to the homomorphism of groups
\[\bZ\times \cM(S;\p)\incl \cM(S;\p)\times \cM(S;\p) \map \cM(S;\p).\] 
and it therefore induces the $S^1$ action $\Phi$ on the quotients. In particular, the $S^1$-action on $|\Fatb|$ 
moves the leaf-edge around the boundary cycle of the metric bordered fat graph.

Fix a bordered fat graph $\Gb$. Assume the leaf-edge of $\Gb$ is attached at a vertex $z$
and that its boundary cycle is $( L \ol{L} A_1 \cdots A_n )$.  Let $x_j$ be the vertex joining $A_j$
and $A_{j+1}$ and let $I_j$ denote the part of the boundary cycle of $\Gb$ lying between the middle of the edge $A_j$ and the middle of the
edge $A_{j+1}$. 

A point in $|\Fatb/\Gb|$ corresponds to a metric fat graph with a morphism of graph to $\Gb$. 
For each such $(f:\tilGb\smap \Gb)$, we identify $S^1$ with the boundary cycle of $\tilGb$ and we
let $I_j(\tilGb)$ denote the part of $S^1$  lying above the segment $I_j$.
For any $\theta\in I_j(\tilGb)$, the circle action attaches the new leaf-edge of $\tilGb$ to a
point above the segment $A_j,A_{j+1}$ in $\Gb$.  
In particular,  we get a map to $\Phi_j(\Gb)$. Hence we have built a map lifting of the
circle action to 
\[\tilde{\Phi} : \big\{(\theta,f:\tilGb\smap\Gb) \tab \theta\in I_j(\tilGb)\big\} \map \Fatb/\Phi_j(\Gb).\]

Other than $x_i$ and $z$, the vertices of $\Gb$ corresponds to a vertex of $\Phi_i(\Gb)$ of the same
valence. If $z$ is trivalent, the lifting $\tilde{\Phi}$ identifies
\[I_i\times \Fat/\Gb \cong \big(I_i\times \cT^{\val(x_i)}\big) \times
\prod_{v\neq x_i,z} \cT^{val(v)}\times \cT^{\val(z)} \maplu{\cong} \Fat/\Phi_i(\Gb)\]
where $I_i$ determines the placement of the new leaf-edge from the
middle of $A_{i-1}$ to the middle of $A_{i}$. Hence, up to signs,
\[\Phi_*([S^1] \tens [\Gb]) = \sum [\Phi_i(\Gb)].\]
If $z$ is not trivalent, the same reasonning shows that the circle action sends the cell at $\Gb$ to
a union of cells obtained by moving the leaf from $z$ to the $x_j$'s. However, in this case, we
do not remove $z$. In particular,
\[\Fat/(\Phi_i(\Gb))\cong \cT^{\val(x_i)+1}\times\cT^{\val(z)-1}\times \prod_{v\neq x_i,z}\cT^{\val(v)}\]
which has the same dimension as $\Fatb/\Gb$. Hence $I_j\times\Fatb/[\Gb]$ is sent to a
cell of lower dimension and its image is degenerate.

We are left to argue our choice of orientation for $\Phi_i(\Gb)$.
For any half-edge $A$ of a bordered fat graph $\tilGb$, let
$\Phi_{A}(\tilGb)$ denote the bordered fat graph obtained by moving
the leaf to the middle of $A$. Hence, the edge $A$ is split into two edges $A_1$ and $A_2$ which are
joined by a trivalent vertex $x$.
Fix a simplex
\[\delta =(\G^b_0\maplu{f_0}\cdots \map \G^b_p=\Gb)\]
and let the boundary component of $\G^b_0$ be $(L\ol{L} D_1\cdots D_m)$.
The circle $\gamma$ at $\G^b_0$ can be subdivided into 
\[ \sum_i (D_i \smap v_i) - \sum_i  (D_i\smap v_{i-1}).\]
where $D_i\smap v_i$ represents the simplex form the middle of $D_i$ to the vertex $v_i$. 
Our choice of orientation is obtained by pushing each simplex $(D_i \smap v_i)\tens \delta$
through to  the following simplex of $[\Phi_{j+1}(\G^b)]$.
\[\nu^+_k = \Phi_{D_i}(\G^b_0)\smap \cdots\smap\Phi_{f_k(D_i)}(\G^b_k)\maplu{D_i^2} \Phi_{f_k(v_i)}(\G^b_k) \smap\cdots
\smap \Phi_{j+1}(\G^b) \]
\end{proof}

\begin{remark}
A similar construction gives a chain complex for the unbordered mapping class group of a surface
$S_{g,n}$ with more than one boundary components. If we consider the mapping class group of $S$ that
fixes only one boundary component then section \ref{sec:CW} gives an integral graph complex to compute its
homology. Again this partially-bordered mapping class group is a $\bZ$ extension of $\cM(S)$. 
\end{remark}

\subsection{Computations for the homology of punctured mapping class groups.}

\begin{proposition}
The homology groups of the punctured mapping class group of the surface
$S_{1,1}$ are
\[H_*(\cM(S_{1,1});\bZ) \cong \begin{cases}
\bZ&*=0\\
\bZ/12 &*=2k+1\\
0 &*=2k+2.
\end{cases}\]
\end{proposition}
\begin{proof}
We have by proposition \ref{prop:S11} that
\[H_*(\cM(S;\p)) = \bZ \qquad * =0,1\]
and is zero for $*\geq2$. Hence the $E^2$-term of the spectral sequence of corollary \ref{cor:ssCPi}
is
\[\begin{array}{ccccccccc}
\bZ& 0 &\bZ& 0 &\bZ& 0& \cdots\\
\bZ& 0 &\bZ& 0 &\bZ& 0& \cdots
\end{array}\]
Consider the cell  $[\G^b_0]$ of figure \ref{fig:objects}. The six bordered fat
graphs $\Phi_i(\G^b_0)$ are isomorphic to $[\G^b_1]$. 
\[\zeta([\G^b_0]) = \sum \Phi_i([\G^b_0]) =\sum \pm [\G^b_1]\]
Since the induced orientations are compatible, we have
\[\zeta([\G^b_0]) = \pm 6 [\G^b_1]\sim \pm 12[\tilG^b_1]\]
where $[\tilG^b_1]$ is a generator of $H_1(\cM(S;\p))$. Hence the $E^3=E^\infty$-term of the spectral sequence is
\[\begin{array}{ccccccccc}
\bZ/12& 0 &\bZ/12& 0 &\bZ/12& 0& \cdots\\
\bZ/12& 0 &0& 0 &0& 0& \cdots
\end{array}\]
and the result follows.
\end{proof}

\begin{proposition}\label{prop:M(S21)}
The homology groups of the punctured mapping class group of the surface
$S_{2,1}$ are
\[H_*(\cM(S_{2,1});\bZ) \cong \begin{cases}
 \bZ &*=0\\
\bZ/10 &*=1\\
\bZ/2\oplus \bZ&*=2\\
\bZ/2\oplus \bZ/120 \oplus \bZ/10&*=2k+3\\
\bZ/2\oplus\bZ/6 &*=2k+4
\end{cases}\]
\end{proposition}
\begin{proof}
Using proposition \ref{prop:s12}, we know that the spectral sequence of corollary \ref{cor:ssCPi}
starts with
\[\begin{array}{ccccccc}
\bZ/6 & 0 & \bZ/6& 0 &\bZ/6 &\dots\\
\bZ/2\oplus\bZ &0 & \bZ/2\oplus \bZ &0&\bZ/2\oplus \bZ & \dots\\
\bZ/2&0&\bZ/2 &0 &\bZ/2 &\dots\\
\bZ/10 & 0&\bZ/10 & 0 &\bZ/10 &\dots\\
\bZ&0&\bZ&0 &\bZ&\dots
\end{array}\]
Using a computer algebra, we find that the circle action is trivial on the homology, which means that $d_2$ is zero.
In fact, using the formula for $d_4$ obtained in the remark~\ref{rmk:ss}, we prove that
the only non-trivial derivatives are 
\[d_4 : E^4_{2n,2q+4}=\bZ \map E^4_{2n-2,2q}=\bZ/2\oplus\bZ\]
which send  $1$ to $(0,120)$. Hence $E^\infty=E^5$ is
\[\begin{array}{ccccccccc}
\bZ/6 & 0 & \bZ/6& 0 &\bZ/6 & 0 &\bZ/6 &\dots\\
\bZ/2\oplus\bZ/120 &0 & \bZ/2\oplus \bZ/120 &0&\bZ/2\oplus \bZ/120 &0 & \bZ/2\oplus \bZ/120 & \dots\\
\bZ/2&0&\bZ/2 &0 &\bZ/2 &0 &\bZ/2 &\dots\\
\bZ/10 & 0&\bZ/10 & 0 &\bZ/10 & 0 &\bZ/10 &\dots\\
\bZ&0&\bZ&0 &0&0&0&\dots
\end{array}\]
This gives the result for $*=0,1,2$ and short-exact sequences
\begin{alignat}{10}
 0 \map \bZ/2\oplus \bZ/120 \map &H_{2k+3}(\cM(S)) \map \bZ/10\map 0 \label{eq:3}\\
0\map \bZ/6 \map &H_{2k+4}(\cM(S))\map \bZ/2 \map 0.\label{eq:4}
\end{alignat}
These extensions are only ambiguous modulo 2 and 5. 

To resolve these extension problems, we consider our spectral sequence with coefficients in $\bZ/2$ and $\bZ/5$. In both
these cases, all derivatives are zero and we get the following $E_\infty$-term.
\[\begin{array}{ccccccccccc}
\bZ/2 & 0 & \bZ/2 & \dots &\qquad\qquad& 0 & 0 &0 &\dots\\
\bZ/2\oplus\bZ/2 & 0 &\bZ/2\oplus\bZ/2&\dots&&0& 0 &0&\dots\\
\bZ/2\oplus\bZ/2\oplus\bZ/2&0&\bZ/2\oplus\bZ/2\oplus\bZ/2 & \dots && \bZ/5&0&\bZ/5&\dots\\
\bZ/2\oplus\bZ/2&0&\bZ/2\oplus\bZ/2&\dots&&\bZ/5&0&\bZ/5&\dots\\
\bZ/2 & 0 &\bZ/2 &\dots && \bZ/5&0&\bZ/5&\dots\\
\bZ/2 & 0 &\bZ/2 &\dots&& \bZ/5 & 0 &\bZ/5& \dots
\end{array}\]
In particular, the second row of the following diagram is an inclusion.
\[ \xymatrix{
\bZ/2\oplus\bZ/2 =H_{3}(\cM(S;\p))\tens \bZ/2 \ar@{^(->}[d]\ar[r]& H_{3}(\cM(S))\tens\bZ/2\ar@{^(->}[d]\\
H_{3}(\cM(S;\p);\bZ/2)\ar@{^(->}[r] &H_{3}(\cM(S);\bZ/2)
}\]
The first is therefore also an inclusion. Using this, tensoring \eqref{eq:3} by $\bZ/2$ gives an exact sequence
\[ 0 \map (\bZ/2)^{\oplus 2} \map H_{3}(\cM(S))\tens \bZ/2 \map \bZ/10\tens \bZ/2 \map 0.\]
In particular,
\[H_{3}(\cM(S))\tens\bZ/2 \cong  (\bZ/2)^{\oplus 3}.\]
Similarly, we have 
\begin{eqnarray*}
H_{3}(\cM(S))\tens \bZ/5 &\cong& (\bZ/5)^{\oplus2}\\
 H_{4}(\cM(S))\tens \bZ/2 &\cong& (\bZ/2)^{\oplus 2}.
\end{eqnarray*}
This proves the claim for $H_3$ and $H_4$.
Since capping with the Euler class $\tau\in H^{2}(\cM(S))$ gives isomorphisms
\[H_{2k+5}(\cM(S))\map H_{2k+3}(\cM(S))\qquad H_{2k+6}(\cM(S))\map H_{2k+4}(\cM(S)),\]
the claim follows completely.
\end{proof}

\subsection{Finite subgroups, symmetries of fat graphs and the homology of $\cM(S)$.}\label{sec:finite}

We shall now illustrate how to use the symmetries of certain (unbordered) fat graphs to find generators for the homology groups of $H_*(\cM(S))$ computed in the previous section.
Fix any finite subgroup $G$ of the mapping class group $\cM(S)$.

\begin{proposition}\label{prop:z/n}
There exists a group isomorphisms $\bZ/n\cong G$ such that the inclusion $f: \bZ/n\smap \cM(S)$ lifts to a map $\tilde{f}:\bZ\smap \cM(S;\p)$ in such a way that the following  diagram commute.
\[\xymatrix{
\bZ\ar[d]^{=}\ar[r]^{n}&\bZ\ar[r]^{1}\ar[d]^{\tilde{f}}&\bZ/n\ar[d]^f\\
\bZ\ar[r]&\cM(S;\p)\ar[r]&\cM(S)
}\]
\end{proposition}
\begin{proof}
For any subgroup $H$ of $\cM(S)$, let $\widetilde{H}$ denote the pullback of the bordered mapping class group over $H$. 
$\widetilde{H}$ is a $\bZ$-extension of $H$ and a subgroup of the torsion-free $\cM(S;\p)$.
Since 
\[H^2(H;\bZ) \equiv [BH,\CPi]\equiv Hom(H,S^1)\]
the extension $\widetilde{H}\smap H$ is characterized by an group homomorphism $\alpha_H:H\smap S^1$. 

We claim that for any finite group $G$, $\alpha_G$ is a monomorphism. Consider $K=Ker(\alpha_G)$. The homomorphism $\alpha_K=\alpha_G|_K =0$ and 
the extension $\widetilde{K}$ is trivial.
In particular, we get
\[K\times\bZ \cong \widetilde{K} \subset \cM(S;\p).\]
Since $\cM(S;\p)$ is torsion free and $K$ is a finite subgroup, we have that $K$ is trivial. Hence $\alpha_G$ is injective and $G$ is a finite subgroup of $S^1$ and cyclic.
Since the extension above $G\cong \bZ/n$ is not trivial, we have
\[ \bZ\map \bZ \map G.\]
We chose the isomorphism $\bZ/n\cong G$ to get the usual $\bZ\smap\bZ/n$.
\end{proof}

For any unbordered fat graph $\G$ in $\Fat$, consider the subgroup
\[\Aut(\G)\subset \pi_1(\Fat)\cong \cM(S).\]
We have an inclusion $\Aut(\G) \incl S^1$
which sends an automorphism $\varphi$ to the distance between $A$ and $\varphi(A)$ along the boundary cycle.
In particular, $\Aut(\G)$ is generated by the automorphism $\varphi_0$ which minimizes this distance.

\begin{corollary}\label{cor:autg}
If $\Gb(X)$ is the bordered fat graph obtained from $\G$ by adding a leaf along the half-edge $X$ then the inclusion $\Aut(\G)\smap\cM(S)$ induces the homomorphism
\[e_1\tens u^k \maps \alpha\tens u^k\]
where $\alpha$ is the path of $\Fatb$ which moves the leaf of  $[\Gb(A)]$ along the boundary cycle until it lies on $\varphi_0(A)$.
\end{corollary}

\begin{proof}
Let $\cB_\bZ$ denote the cellular chain complex
\[\bZ \imaplu{0} \bZ \imaplu{0} 0 \dots\]
for the circle. Let the chain complex $D_{n}=\cB_\bZ\tens\bZ[u]$ have derivative
\[d_n(e_1\tens u^k) = 0\qquad 
d_n(e_0\tens u^k) = \begin{cases}
n(e_1\tens u^{k-1}) &k\geq0\\
0&k=0.
\end{cases}\]
This is the chain complex of section \ref{sec:DX} associated to the action of the circle on itself by multiplication by $z^n$. It therefore computes the $S^1$-equivariant homology of this action which gives the homology of $\bZ/n$.

Using the naturality of this construction, the generator $e_1\tens u^{k}$ of $H_{2k+1}(\bZ/n)$ is sent to the cycle
$\tilde{f}(1)\tens u^k$ where $\tilde{f}$ is the lifting promised by proposition \ref{prop:z/n}. 
The image of each morphism in $\alpha$ gives an isomorphism of $\Fat$ between two fat graphs isomorphic to $\G$. Their composition  is homotopic to an automorphism of $\G$ which moves the edge $A$ which supported the leaf in $\G(A)$ to the edge $\varphi_0(A)$ which supported the leaf in $\G(\varphi_0(A))$.
It is left to show that the element $\alpha$ is the lifting $\tilvarphi$ promised by proposition \ref{prop:z/n}.  This follows easily from the fact that $\alpha^n$ moves the leaf once around the boundary component.
\end{proof}

Consider the two group homomorphism
\[\bZ/6 \maplu{\kappa_6} \cM(S_{1,1}) \imaplu{\kappa_4}  \bZ/4 \]
obtained by including the automorphism groups of the unbordered fat graphs $\G_6$
and $\G_4$ with boundary cycle
\[\om_1=(ABC\ol{ABC})\qquad \om_2=(EF\ol{EF}).\]
Note that these homomorphism are obtained by building the surface  $S_{1,1}$ out of  a onced-punctured hexagon and a onced-punctured square.

\begin{corollary}
The homology group $H_{2k+1}(\cM(S_{1,1}))$ is generated by 
\[\epsilon_k = \kappa_4(e_{4,k}) - \kappa_6(e_{6,k})\]
where $e_{i,k}$ is a generator of $H_{2k+1}(\bZ/i;\bZ)$.
\end{corollary}
\begin{proof}
Using both proposition \ref{prop:DX} and corollary \ref{cor:autg}, we have that
\[\epsilon_k= (\psi_4-\psi_6)\tens u^k.\]
Since capping with $\tau$, which induces an isomorphism, sends $\epsilon_k$ to $\epsilon_{k-1}$, it suffices to show that $\psi_4-\psi_6$ generate $H_1(\cM(S;\p))$ modulo $12$.

The lifting $\psi_6$ corresponds to the cell $[\G^b_1]$ of figure \ref{fig:s11} of page \pageref{fig:s11}. Since $\G_4$ has dimension one, the lifting $\psi_4$ is not in the graph complex. Using the notation of figure \ref{fig:s11}, the diagram
\[\xymatrix{
\tilG^b_1\ar[rr]^{N}&&\G^b_2&&\ar[ll]_M \tilG^b_2\\
\G^b_0\ar[u]^H \ar[urr]^{E,H} \ar[rr]_{E}&&\G^b_1\ar[u]^{K}&&\ar[ll]^F \ar[ull]_{F,G} \ar[u]_{G} \G_0^b
}\]
commutes. Its top row is $\psi_4$; its bottom one is $[\G^b_1]$ and its vertical part gives $[\tilG^b_1]$. We therefore get that 
\[\psi_4-\psi_6 = [\tilG^b_1] + [\G^b_1] - [\G^b_1] = [\tilG^b_1]\]
 which generated $H_1(\cM(S;\p))$. 
\end{proof}

We will now attack the more interesting case of the homology of the unbordered mapping class group of the surface $S_{2,1}$.
Consider the fat graphs $\G_{10}$ and $\G_2$ which have the following boundary cycles.
\begin{equation}\label{eqn:g2}
(ABCDE\ol{ABCDE}) \qquad
(GHIGJHIJFKMNKPMNPF).
\end{equation}
The second graph can be obtained by removing the leaf of the first graph in figure \ref{fig:gen HS12}.
These give group homomorphisms
\[\bZ/10\maplu{\kappa_{10}} \cM(S_{2,1})\imaplu{\kappa_2}\bZ/2.\]

\begin{proposition}Ê\label{prop:z10}
The element
\[ \beta_{k} = \kappa_{10}(e_{10,k}) -\kappa_{2}(e_{2,k})\]
generates the $\bZ/10\subset H_{2k+1}(\cM(S))$ of highest filtration.
\end{proposition}
\begin{proof}
Again it suffices to show that $\beta_1$ generates $H_1(\cM(S))$. As in corollary \ref{cor:autg}, we lift the generator of $\Aut(\G_2)$ to the cycle
\[\begin{split}
\big[(LLGHIGJHIJFKMNKPMNPF)\big] &+  \big[(LLHIGJHIJFKMNKPMNPFG)\big]\\
+\dots &+\big[(LLFKMNKPMNPFGHIGJHIJ)\big].
\end{split}\]
The lifting $\alpha$ of $\bZ/10$ to a $1$-cycle of $\Fatb$ does not land in the graph complex. We first expand $\G_{10}$ to a trivalent graph $\tilG_{10}$ with boundary cycle
\[(AWXBCXDZEAZYBCYDWE).\]
For any fat graph $\G$, let $\G(X)$ be the bordered fat graph obtained by adding a leaf on the half-edge $X$. By pushing the leaf movement to $\tilG_{10}$, we get that $\alpha$ corresponds to the following chain of $\Fatb$
\[ \begin{split}&
\big[(LLWXBCXDZEAZYBCYDWEA)\big]+ \big[(LLXBCXDZEAZYBCYDWEAW)\big]\\
&+\big[(LLBCXDZEAZYBCYDWEAWX)\big]-\varphi+\psi.
\end{split}\]
where $\varphi$ and $\psi$ are the natural morphisms $[\tilG_{10}(A)\smap \G_{10}(A)]$ and $[\tilG_{10}(B)\smap\G_{10}(B)]$. Using that $\G_{10}(A)$ and $\G_{10}(B)$ are isomorphic, $\varphi$ and $\psi$ lifts to the cell $\Fatb/[\G_{10}(A)]$. We find that the 1-cells of $\Fatb$
\begin{equation}\begin{split}\label{eqn:phipsi}
&\big[(LLA^\prime BCZDEZYA^\prime ABYCXDEXA)\big] \\
&+ \big[(LLA^\prime BCZDEZYA^\prime AWBYCWDEA)\big]\\
&+\big[(LLA^\prime BCZDVEZA^\prime AVWBCWDEA)\big] \\
&+\big[( LLA^\prime VWBCWDEA^\prime AUBCUDVEA)\big]
\end{split}\end{equation}
lift to connect $\varphi$ and $\psi$ in $\Fatb/[\G_{10}(A)]$. Using that the geometric realization of $\Fatb/[\G_{10}(A)]$ is contractible, we get that \eqref{eqn:phipsi} is homotopy equivalent to $\psi-\varphi$. Hence we have pushed $\alpha$ to the graph complex. Using magma, we can now show that $\epsilon_1$ generates $H_1(\cM(S;\p))$.
\end{proof}

Consider the graph $\G$ with boundary
\[\om = (AB{C}D{EA}FC{G}H{DF}IG{B}E{HI}).\]
Its automorphism group is isomorphic to $\bZ/3$

\begin{proposition}
The cycle $\theta_k= \kappa_3(\epsilon_{3,k})$ generates the $\bZ/3$-part of $H_{2k+3}(\cM(S_{2,1}))$.
\end{proposition}
\begin{proof}
Since $H_1(\cM(S))\cong\bZ/10$, $\kappa_3$ induces the zero map between the first homology groups.The lifting $\alpha$ predicted by corrolary \ref{cor:autg} gives 
\[\begin{split}
&\big[(LL ABCDEAFCGHDFIGBEHI)\big] + \big[(LLBCDEAFCGHDFIGBEHIA)\big]\\ 
&+\big[(LLCDEAFCGHDFIGBEHIAB)\big]+\big[(LLDEAFCGHDFIGBEHIABC)\big]\\ 
&+\big[(LLEAFCGHDFIGBEHIABCD)\big]+ \big[(LLAFCGHDFIGBEHIABCDE)\big]
\end{split}\]
in $\GC_1$. Let $\beta\in\GC_2$ be a chain that kills $\alpha$. Since
\[d_{tot}(\beta\tens u) = \alpha\tens u + \zeta(\beta)\tens 1,\]
$\alpha\tens u$ is homologous to the image in $\Fat$ of $\zeta(\beta)$. Using magma, we find $\beta$ and show that $\zeta(\beta)$ is
homologous to $40$-times the generator of 
\[\bZ\subset H_3(\cM(S;\p))\map \bZ/120\subset H_3(\cM(S))\]
and hence maps to the $\bZ/3$ part of $\bZ/120$. 
\end{proof}

The automorphism group of the fat graph $\G_8$ with boundary cycle
\[ \om=(ABCDABCD).\]
is isomorphic to $\bZ/8$ with generator $\varphi_8$. We consider
\[\bZ/8\maplu{\kappa_8} \cM(S).\]

\begin{proposition}
The cycles
\[\theta_k=\kappa_{8}(\epsilon_{8,k})-\kappa_{2}(\epsilon_{2,k})\]
generate the $\bZ/8$-part of $H_{2k+3}(\cM(S))$. Here $\kappa_2$ is the homomorphism considered in proposition \ref{prop:z10}.
\end{proposition}
\begin{proof}
The lifting of $\epsilon_{2,k}$ was found in the proof of proposition \ref{prop:z10}. As before we find that the lifting promised by corollary \ref{cor:autg} splits into a movement of leaves on a trivalent graph above $\G_{8}$
and a path in one cell of $\Fatb$.
Using magma, we show that $\theta_1$ is trivial. 

Using the morphism
\[\tilG=(AEFGBCGHDAHFIBCIED)\map (ABCDABCD),\]
we can show that $(\varphi_8)^4$  is conjugate to the only non-trivial automorphism of $\tilG$. Since $\tilG$ is trivalent, corollary \ref{cor:autg} gives us a lifting $\tilde{\alpha}$ in $\GC_*$.  Using magma and the reasonning of the previous proposition, we show that $\tilde{\alpha}\tens u$ generates
\[\bZ/2\subset\bZ/8\subset H_{2k+3}(\cM(S)).\]
\end{proof}

\section{Infinite loop operations and the bordered graph complex}
\label{sec:ops}

In \cite{Tillmann}, Tillmann showed that the stable homology of the bordered mapping class group was the homology
of an infinite loop space. In particular, it comes equipped with Araki-Kudo-Dyer-Lashoff operations.
The first of these operations comes from Miller's algebra over the little-disk operad.

\subsection{Pair of pants product.}\label{sec:Miller}

By gluing a chosen boundary component of two surfaces $S_{g,k}$ and $S_{h,j}$ to the first and second boundary components of a
pair of pants, we get a homomorphism
\[ \cM(S_{g,k};\p)\times \cM(S_{h,j};\p)\map \cM(S_{g+h,k+j-1};\p).\]
We will construct a combinatorial version of this structure.

Take two bordered fat graphs $\Gb$ and $\tilGb$ with chosen boundary cycles
\[\om_1=(L_1\ol{L}_1 A_1\ldots A_n)\qquad\qquad \tilom_1=(L_2\ol{L}_2 B_1\ldots B_m). \]
Denote by $\Gb\#\tilGb$ the bordered fat graph obtained by identifying
the two leaf-vertices to a new vertex $x$. We attach a new leaf-edge
$L_\new$ to $x$ so that the first boundary cycle of this new bordered
fat graphs becomes   
\[\om = (\ol{L}_1 A_1\ldots A_n\, L_1\, L_{\new}\,\ol{L}_{\new} \ol{L}_2 \, B_1\ldots B_m\, L_2).\]
Fix orientations
\begin{eqnarray*}
o_{\Gb} &=&(v_1\wedge \ldots \wedge v_r\wedge l_1) \tens(a_1\wedge\ldots \wedge a_n\wedge L_1\wedge \ol{L}_1)\\
o_{\tilGb}&=&(w_1\wedge\ldots \wedge w_s\wedge l_2)\tens (b_1\wedge\ldots \wedge b_m\wedge L_2\wedge \ol{L}_2).
\end{eqnarray*}
We construct the following orientation on $\Gb\#\tilGb$.  
\[\begin{split} 
o_{\Gb}&\# o_{\tilGb}= (w_1\wedge\ldots\wedge w_s)\wedge x \wedge
(v_1\wedge \ldots \wedge v_r)\wedge l\\ 
&\tens(b_1\wedge\ldots \wedge b_m)\wedge L_2 \wedge \ol{L}_2 \wedge
(a_1\wedge\ldots \wedge a_n)\wedge L_1 \wedge \ol{L}_1 \wedge 
L_\new\wedge \ol{L}_\new 
\end{split}\]

\begin{proposition} \label{prop:K}\label{prop:pants}
The map
\begin{eqnarray*}
\mu_*\ :\ \GC_p \tens \GC_q &\map& \GC_{p+q}\\ 
{[\Gb,o_\G]}\tens [\tilGb,o_\tilG] &\map& [\Gb\#\tilGb,o_\G\#o_\tilG]
\end{eqnarray*}
induces the pair of pants gluing on the homology of the mapping class groups.
\end{proposition}
\begin{proof}
The functor
\[\mu( \Gb , \tilGb)= \Gb\# \tilGb\]
glues two bordered fat graphs along the bordered fat graph
of $\Fat^b_{0,3}$ with boundary cycles
\[\om=(L_0\,L_0\,C\,\ol{B}\,\ol{C}\, D\,\ol{A}\,\ol{D})(L_1\,L_1\,A)(L_2\,L_2\,B).\]
By section \ref{sub:complex}, the cell at $\Gb\#\tilGb$ is the image of the category
\[ \Fatb/[\Gb\#\tilGb] \cong \prod_{v \in \Gb\#\tilGb} \cT^{val(v)}.\]
Every interior vertex of $\Gb$ and of $\tilGb$ corresponds naturally to
a vertex of $\Gb\#\tilGb$ of the same valence. In fact, $\mu$ identifies 
\[\begin{split}
\Fatb/[\Gb] \times\Fatb/[\tilGb] &\cong  \prod_{v\in\Gb} \cT^{val(v)}\times \prod_{w\in \tilGb} \cT^{val(w)}\\
&\maplu{\cong}  \left(\prod_{v\neq x} \cT^{val(v)}\right)\times
\cT^{val(x)} \cong \Fatb/[\Gb\#\tilGb]. 
\end{split} \]
Hence $\mu$ induces a cellular map on $\Fatb$ which gives the pair of pants gluing.

It now suffices to verify our choice of orientation. Fix simplices
\[(\delta=\G^b_0\smap\cdots \smap \G^b_p)\qquad (\lambda=\tilG^b_0\smap\cdots\smap\tilG^b_q).\]
Fix also a $(p,q)$-shuffle $\kappa=(\kappa_1,\kappa_2)$, ie two
increasing maps 
\[\kappa_1:\{0\dots p+q\} \map \{0\dots p\}\qquad \kappa_2:\{0\dots
p+q\}\map \{0\dots q\}\]
so that $\kappa_1+\kappa_2=Id$. The image of the chain
\mbox{$\delta\tens\lambda$} hits each simplex 
\[\nu_\kappa = (\G^b_{\kappa_1(0)} \#  \tilG^b_{\kappa_2(0)})\map\cdots\map (\G^b_{\kappa_1(p+q)} \#
\tilG^b_{\kappa_2(p+q)})\]
with sign $(-1)^{|\kappa|}$ coming from the Eilenberg-Zilber map. It
therefore suffices to show that 
\begin{equation}\label{eqn:compatible}
 (-1)^{|\kappa|} \epsilon(\delta, o_\Gb)\ \epsilon(\lambda, o_\tilGb)= \epsilon(\nu_\kappa, o_\Gb\# o_{\tilGb})
\end{equation}
where, as in lemma \ref{lem:orient pcells}, $\epsilon(\alpha,o)$ equals $\pm 1$ depending on whether
the orientation $o$ agrees or disagrees with the orientation induced by the natural orientation on
the trivalent graph in $\alpha$. By remark \ref{rmk:chain}, inverting two consecutive morphisms
changes the right hand side of \eqref{eqn:compatible} by $-1$. However this also changes the sign of
$\kappa$. Hence, we only need to verify \eqref{eqn:compatible} for the simplex
\[\nu = ({\G^b_0}\# \tilG^b_0) \smap\cdots\smap ({\G^b_p}\# \tilG^b_0) \smap\cdots\smap ({\G^b_p}\#\tilG^b_q)\]
which is straightforward.
\end{proof}

\begin{proposition}
The pair of pants product is part of an $A_{\infty}$ structure on
$X=|\Fat|_+$. In particular, this product is homotopy associative in a natural way.
\end{proposition}
\begin{proof}
Using theorem 5 of \cite{Stasheff}, we only need to show that they are compatible functors 
\[\mu_k:\cT^{k+1} \times \big(\Fat^b\big)^k\map \Fat^b\]
where $\cT^k$ is the category of plannar trees used in section \ref{sec:CW}.
We define
\[\mu_k\big([T_k],\G_1,\dots,\G_k\big)\]
to be the bordered fat graph obtained by identifying the first leaf-edge of the
bordered fat graph $\G_j$ to the \xth{j} leaf-edge of $[T_k]$. The zero$^{th}$ leaf
of $[T_k]$ becomes the first leaf of the new bordered fat graph.
\end{proof}

\begin{example}
The pair of pants gluing
\[\cM(S_{1,1};\p))\times \cM(S_{1,1};\p)\map \cM(S_{2,1};\p).\]
send the generator of figure \ref{fig:genS11} to the generators of 
\[H_0(\cM(S_{2,1};\p))\cong\bZ \quad H_1(\cM(S_{2,1};\p))\cong \bZ/10\quad
H_2(\cM(S_{2,1};\p))\cong\bZ/2\]
as shown in figure \ref{fig:genS21}. 
\end{example}

\begin{example}
The pair of pants gluing
\[\cM(S_{1,1};\p))\times \cM(S_{0,2};\p)\map \cM(S_{1,2};\p).\]
sends the generator of figure \ref{fig:genS11} and \ref{fig:finb12} to the first four generators of figure \ref{fig:gen HS12}. 
\end{example}

\subsection{First Dyer-Lashoff-Kudo-Araki and Browder operations.}

The Pontrjagin product on a loop space is homotopy associative. The Dyer-Lashoff-Kudo-Araki measure how far it is from being commutative at the chain level. The Browder operations of an $n$-fold loop space $X$ give obtructions to $X$ being in fact an $n+1$-fold loop space.

Let $\cC(*)$ denote the little 2-disks operad. An element of $\cC(j)$ consists of $j$ disks within the unital
disks. Miller in \cite{Miller} pointed out the existence of maps 
\[\Omega : \cC(j) \times \cM(S_{g_1,1};\p)\times\ldots\times \cM(S_{g_j,1};\p)\map \cM(S_{(\sum g_j),1};\p)\]
obtained by gluing the surface $S_{g_i,1}$ to the \ith disk. He then used these maps to show that 
\[X= \Du_{g} B\cM(S_{g,1};\p).\]
is a $\cC(2)$-algebra and that its group completion $X^+$ is a two-fold loop space. 
This structure give us the first Kudo-Araki operation and the first Browder operation.

\begin{remark}
Since the stable mapping class group has the homology of an infinite loop space,
any Browder operation on $X$ is stably trivial. However Fiedorowicz and Song
have shown in \cite{Fied_Song} that 
\[\Psi_1 : H_0(\cM(S_{1,1};\p))\tens H_0(\cM(S_{1,1};\p))\map H_1(\cM(S_{2,1};\p))\]
is non-trivial before stabilization.
\end{remark}

We define a chain map $q_1$ which sends $\sum  a_j [\G^b_j,o_j]$ to
\[\sum_{j,k} a_ja_k\Bigg(\mu\Big([\G^b_j,o_j],\zeta\big([\G^b_k,o_k]\big)\Big)- \sum_{i\in
  \G^b_j} \Phi_i\Big(\mu\big([\G^b_j,o_j],[\G^b_k,o_k]\big)\Big) \Bigg)\]
where the last sum covers only the vertices of $\G^b_j$. Here $\mu$ is the pair of pants product, $\zeta$ is the circle action and $\Phi_i$ was used in section \ref{sec:DX} to define $\zeta$

\begin{proposition}\label{prop:kudo}
The map $q_1$ induces the Kudo-Araki operation 
\[Q_1:H_{k}(\cM(S_{g,1};\p);R)\map H_{2k+1}(\cM(S_{2g,1};R))\qquad \]
where $R$ is either $\bZ/2$ or $\bZ$ depending on the parity of $k$.
\end{proposition}

\begin{figure}
\begin{center}
\mbox{
\subfigure[{$\Gb \#_{(t_0,t_1,t_2)} \Gb$} ]{\quad\epsfig{file=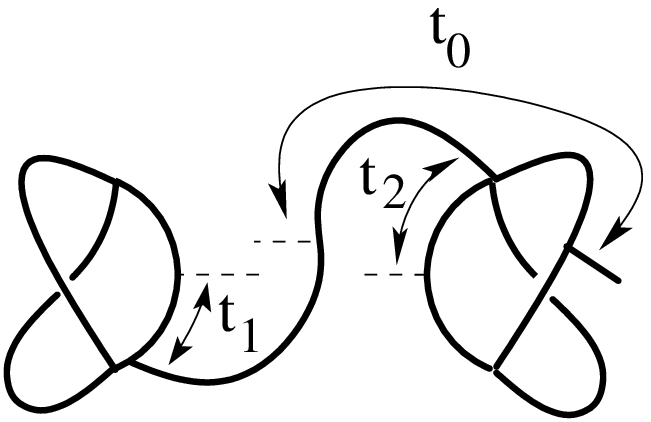, height=55pt}\quad}
\subfigure[{$S^1\smap\cC(2)\smap |\Fat^b|$}]{\quad\epsfig{file=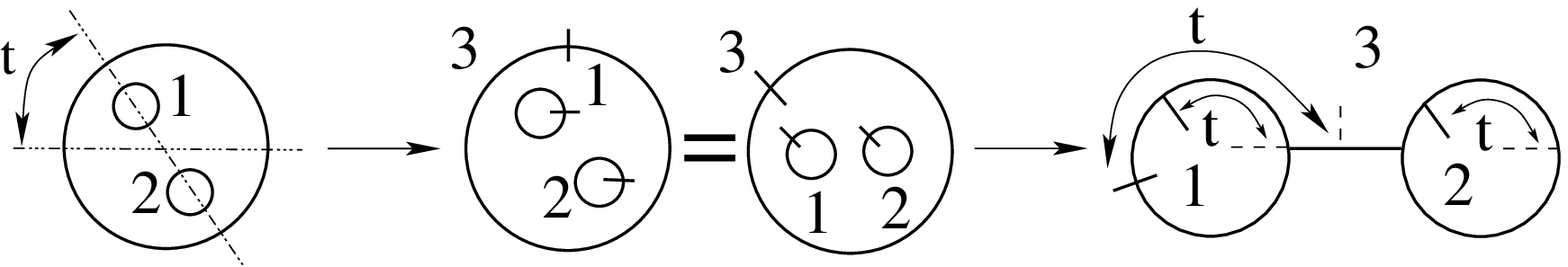, width=220pt, height=40pt}\quad}}
\caption{A graph version of Miller's construction}\label{fig:Miller}
\end{center}
\end{figure}

\begin{proof}
The space $|\Fat^b_{0,3}|$ is homotopy equivalent to $(S^1)^3$ each circle
representing the movement of one of the leaf around its boundary cycle.
We define a parameterized version of the pair of pants gluing $\mu$ as follows.
\begin{eqnarray*}
\xi\ :\  |\Fat^b_{0,3}| \times |\Fat^b_{g,1}|\times |\Fat^b_{h,1}| &\map& |\Fat^b_{g+h,1}|\\
((t_0,t_1,t_2), \Gb,\tilGb) &\maps& \Big(\Gb \underset{t_0,t_1, t_2}\# \tilGb\Big).
\end{eqnarray*}
The new metric bordered fat graph, which illustrated in figure \ref{fig:Miller},
is obtained by removing the leaf-edges of $\Gb$ and $\tilGb$ and by adding a new
edge  $K$ of length one which is attached $t_1$-th of the way around the
boundary cycle of $\Gb$ and $t_2$-th of the way around the boundary
cycle of $\tilGb$. Finally a leaf is attached $t_0$-th of the way
around the new boundary whose starting point is middle of the edge $K$.

As in figure \ref{fig:Miller}, we construct a $\bZ/2$-equivariant map,
\[f\ :\ S^1\times X\times X \map \cC(2)\times X\times X \map |\Fat^b_{0,3}|\times X\times X \maplu{\xi} X.\]
The group $\bZ/2$ acts on these spaces : it exchanges the two
coordinates of $X\times X$, the two disks of $\cC(2)$ and two of the boundary
components of $S_{0,3}$.
Let $e$ be the one-simplex corresponding to the top half of $S^1$. For any $k$-cycle $u$ of X, the chain $e\tens u\tens u$ maps to a cycle of $C_{2k+1}(X;R)$. The Kudo-Araki is defined to be the homology class represented by this cycle.

The map
\[C_1(S^1)\tens C_k(X)\map C_1((S^1)^3)\tens  C_k(X)^{\tens 2} \maplu{\cong} C_1(S^1)^{\tens 3}\tens C_k(X)^{\tens 2}\]
sends $e\tens u$ to 
\[(\pi\tens e\tens \pi\tens u\tens u) +((-e)\tens 0\tens\pi\tens u\tens u) + 
(0\tens 0\tens e\tens u\tens u)\]
Let $\gamma \in C_1(S^1)$ be the entire circle. Using that $f$ is $\bZ/2$-equivariant, 
\[Q_1(u) = \xi(0\times 0\times \gamma\times u\times u) + \xi(-e\times 0
  \times\pi\times u\times u).\]
Using proposition \ref{prop:pants} and lemma \ref{lem:circle}, this gives the result.
\end{proof}

\begin{example}
The Dyer-Lashoff operation 
\[ Q_1 : H_*(\cM(S_{1,1};\p)) \map H_*(\cM(S_{2,1};\p))\]
maps the generator $[\G^b_0]$ for $H_0(\cM(S_{1,1};\p);\bZ)$ to
\[\begin{split}
q_1([\G^b_0])&= \mu_*([\G^b_0],\zeta([\G^b_0])) +\sum_i\Phi_i([\G^b_0\#\G^b_0])\\
&= -6* [\G^b_0\#\G^b_1]+ \sum_i(\Phi_i([\G^b_0\#\G^b_0])) = -12+5=3.
\end{split}\]
To see that the second term gives $5$, notice that it is the lifting of $Aut(\G_2)$ of \eqref{eqn:g2} to $\GC_*$. The result then follows from the proof of proposition \ref{prop:z10}.
Similarly $q_1$ maps $H_1(\cM(S_{1,1};\p);\bZ/2)$ surjectively to the cokernel of the map
\[H_3(\cM(S_{2,1};\p);\bZ)\tens \bZ/2 \incl H_3(\cM(S_{2,1};\p);\bZ/2).\]
\end{example}

We define
\[\psi_1\ :\ \GC_p \tens \GC_q \map\GC_{p+q+1}\]
by sending $[\Gb,o]\tens[\tilGb,\tilde{o}]$ to
\[-\xi\Big(\mu\big([\Gb,o],[\tilGb,\tilde{o}]\big)\Big)
+\mu\Big(\zeta([\Gb,o]),[\tilGb,\tilde{o}]\Big)
+(-1)^p \mu\Big([\Gb,o], \zeta\big([\tilGb,\tilde{o}]\big)\Big).\]

\begin{proposition}
The map $\psi_1$ induces the Browder operation on the homology of the
bordered mapping class groups. 
\end{proposition}

\begin{proof}
Using the notation of the proof of proposition \ref{prop:kudo}, the Browder operation is defined by
\[\psi_1({u\tens v}) = f(\gamma\tens u\tens v)\]
where again $\gamma\in C_1(S^1)$ is the entire circle. Since $\gamma$ is first mapped to 
\begin{eqnarray*}
 (\gamma^{-1} \times\gamma \times (\gamma+\pi)) &\simeq& (\gamma^{-1}\times 0\times \pi) +(0\times
\gamma\times \pi) + (0\times 0 \times \gamma +\pi)\\
&\simeq& (\gamma^{-1}\times 0\times 0) +(0\times \gamma\times 0) + (0\times 0 \times \gamma )
\end{eqnarray*}
in $(S^1)^3$, it suffices to find the image under $\xi_*$ of the chains
\[(\gamma^{-1}\times 0\times 0 \times [\Gb] \times [\tilGb]) \quad (0\times \gamma\times 0 \times
     [\Gb] \times [\tilGb])
\quad (0\times 0\times \gamma \times [\Gb] \times [\tilGb])\]
which is done as in the previous proof.
\end{proof}
\begin{example}
The Browder operation $\Psi_1$ gives
\begin{eqnarray*}
H_0(\cM(S_{1,1};\p))\tens H_0(\cM(S_{1,1};\p)) &\maplu{6}& H_{1}(\cM(S_{2,1};\p))\\
H_0(\cM(S_{1,1};\p))\tens H_1(\cM(S_{1,1};\p)) &\maplu{0}& H_{2}(\cM(S_{2,1};\p))\\
H_1(\cM(S_{1,1};\p))\tens H_1(\cM(S_{1,1};\p)) &\maplu{0}& H_{3}(\cM(S_{2,1};\p)).
\end{eqnarray*}
In particular, this shows that $6\in H_1(\cM(S_{2,1};\p))$ is unstable. 
\end{example}

\subsection{Higher Kudo-Araki-Dyer-Lashof operations.}

Let $S_{g,1}$ be a surface of genus $g$ with a single boundary component. As before, let $\cM_{\infty}$ be the direct limit of the following sequence of groups.
\[\cM(S_{1,1};\p) \maplu{\Phi_1} \cM(S_{2,1};\p) \maplu{\Phi_2} \cM(S_{3,1};\p)\map \dots\]
The infinite loop structure of Tillmann \cite{Tillmann} on the homology of the stable mapping class groups gives rise to Dyer-Lashoff-Araki-Kudo operations
\[\widetilde{Q}_{i,p} : H_n(\cM_{\infty,\delta};\bZ/p) \map H_{pn+i}(\cM_{\infty,\delta};\bZ/p).\]
Before the existence of such an infinite loop structure was known, Cohen and Tillmann constructed homological operation in \cite{Cohen_Tillmann} using simple geometric construction and Harer stability. These operations are the Kudo-Araki-Dyer-Lashof operations associated to the infinite loop structures on $\cM_\infty$. 

In this section, we construct an unstable and combinatorial version of these operations. These will give homomorphisms
\[Q_{i,p}: H_n(\cM(S_{g,1};\p);\bZ/p) \map H_{pn+i}(\cM(S_{pg,1});\bZ/p).\]
The infinite loop space operations are obtained from these by composing with the stable splitting
\[s: H_n(\cM(S))\map H_n(\cM(S;\p))\]
which exists from the proof of corollary \ref{cor:trivialS1}. 

\begin{figure}
\begin{center}
\mbox{\epsfig{file=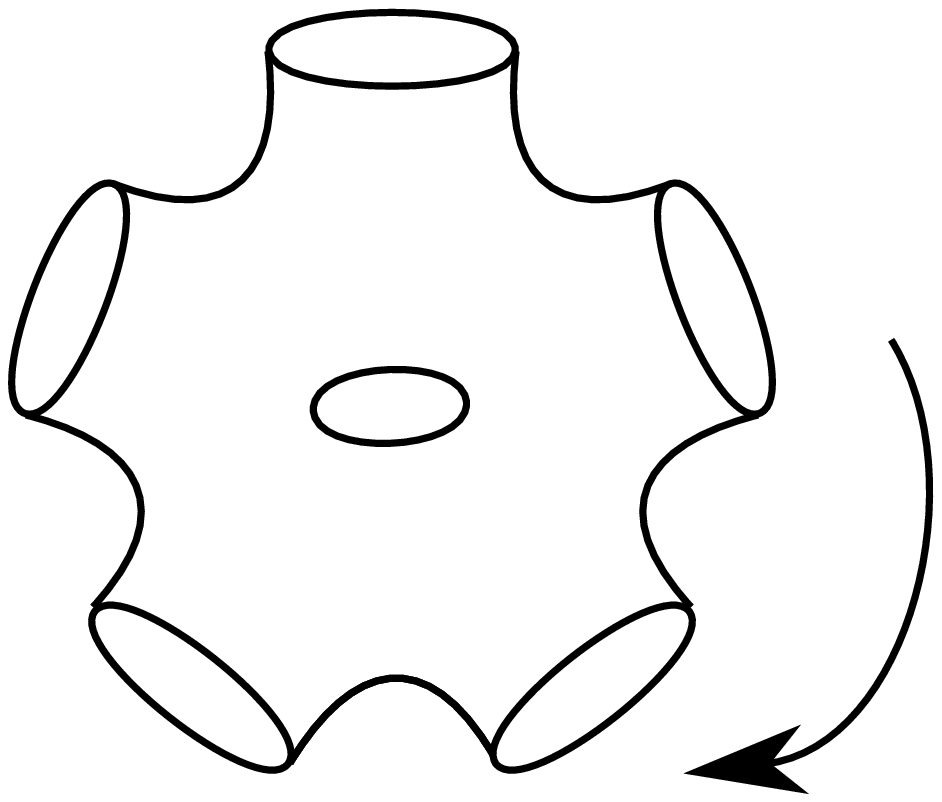, width=75pt}}
\caption{The embedding $\bZ/5\smap \ctilM_{0,5}^1$}\label{fig:S06}
\end{center}
\end{figure}

Fix a prime $p$ and a surface $S_{0,p+1}$ of genus 0 with $p+1$ boundary components. Fix parameterizations $\{\phi_i\}$ of the first $p$ boundary components. Consider the group 
\[\tilm=\pi_0\Diff(S_{0,p+1};\{\phi_i\})\] 
of isotopy classes of diffeomorphisms of $S$ which preserve the last boundary component as a set. These diffeomorphisms are allowed to permute the first $p$ boundary components but they must preserve the parameterizations. 
Using the action of the elements of $\tilm$ on the first $p$ boundary components, we get a wreath product
\[\tilm\int \cM(S_{g,1};\p) = \tilm\ltimes \big(\cM(S_{g,1};\p)^p\big).\]
The rotation illustrated in figure \ref{fig:S06} gives an inclusion $\bZ/p\smap\tilm$. 
By gluing the $p$ surfaces $S_{g,1}$ onto the first $p$ boundary components of $S_{0,p+1}$, we get
\begin{equation}\label{eqn:twistedprod}
\bZ/p\int \cM(S_{g,1};\p)\map \tilm\int  \cM(S_{g,1};\p) \map \cM(S_{pg,1})
\end{equation}
which gives a chain map
\[\theta : C_*(E\bZ/p)\underset{\bZ/p}{\tens}\big(C_*(B\cM(S_{g,1};\p))\big)^{\tens p}
\maplu{\Psi_*} C_*( \cM(S_{pg,1}))
\]
For a cycle $\alpha\in C_{k}(\cM(S_{g,1};\p);\bZ/p)$, we
define
\[Q_i([\alpha])= \theta\big[e_i\tens x^{\tens p}\big] \in
H_{pk+i}(\cM(S_{pg,1};\p);\bZ/p)
\]
where each $e_i$ is one of the generator of the standard resolution of $\bZ/p$ over itself.

Consider the following category $\cC(p,g)$ which is the wreath product of the category $\Fatb$ with one-object category $\cC(p)$ whose morphisms are $\bZ/p$. 
The objects of $\cC(p,g)$ are $p$-tuples of bordered fat graphs. The morphisms of $\cC(p,g)$ are tuples 
\[(i,\varphi_1\dots\varphi_p) : (\G^b_1,\dots \G^b_p)\map (\tilG^b_1,\dots \tilG^b_p)\] 
Here $\varphi_j:\G^b_j\smap \tilG^b_{j+i}$ is a morphism of bordered fat graphs. 
By construction
\[|\cC(p,g)| \simeq B\left( \bZ/p\int \cM(S;\p)\right)\]

Following the contruction of the pair of pants functor $\mu$, we define
\[\mu_p : \cC(p,g)\map \Fat_{pg}^1.\]
A $p$-tuple $(\G^b_1\dots \G^b_p)$ is sent to the unbordered fat graph $\G_\new$
obtained by identifying the leaf-vertices $l_1,\dots, l_p$. 
A morphism of $\cC(p,g)$ induces a morphism of $\Fat$ in the obvious way. The functor $\mu_p$ realizes the homomorphism of groups of \eqref{eqn:twistedprod}.
We therefore get the homomorphism $\theta$ by first using the shuffle map 
\[\begin{split}
C_*(\cE(p))\underset{\bZ/p}\tens \big(C_*(\Fat^b_{g,1})\big)^{\tens p}& \map C_*\left(\cE(p)\underset{\bZ/p}\tens\big(\Fat^b_{g,1}\big)\right) = C_*(\cC(p,g))\\
&\map C_*(\Fat_{pg,1})
\end{split}\]
Here $\cE(p)$ is the universal cover of $\cC(p)$. It has $p$ object and a single morphism between any two of them. 
The chains $e_i\in C_k(\cE(p);\bZ/p))$ which are the image of the generators of
the standard resolution are as follows.
\begin{eqnarray*}
 e_{2k}&=& \sum_{1\leq i_1\dots i_k\leq p-1} (i_1,1,i_2,1\dots , i_k,1) =
(0\smaplu{i_1}i_1\smaplu{1}(i_1+1)\smaplu{i_2}\dots ) \\
 e_{2k+1}&=& \sum_{i_1\dots i_k} (1, i_1, 1, i_2,\dots, i_k,1) =
(0\smaplu{1} 1\smaplu{i_1}(i_1+1)\smaplu{1}(i_1+2)\dots ) \\
\end{eqnarray*}

\begin{example}
Consider the operations
\[Q_{i,2}:H_{0}(\cM(S_{1,1};\p);\bZ/2)\map H_{i}(\cM(S_{2,1});\bZ/2).\]
Let $[\G^b_0]$ be the generator of $H_0(\cM(S_{1,1};\p)$, let $\tilG_2$ be again the fat graph of page \pageref{eqn:g2}) and let $\varphi$ be the generator of $\Aut(\tilG_2)\cong \bZ/2$. By definition 
\[Q_{i,2}([\G^b_0]) = [(\varphi,\varphi,\dots,\varphi)] \in H_{i}(\Fat;\bZ/2)\]
In proposition \ref{prop:z10}, we have shown that for $i$ odd
\[  [(\varphi,\varphi,\dots,\varphi)] =5\in \bZ/10\subset H_*(\cM(S_{2,1};\bZ)\]
By the proof of proposition \ref{prop:M(S21)}, we know that the $E^\infty$ term of the appropriate Lerray-Serre spectral sequence with $\bZ/2$-coefficients is
\[\begin{array}{ccccccccccc}
\bZ/2 & 0 & \bZ/2 &\cdots\\
\bZ/2\oplus\bZ/2 & 0 &\bZ/2\oplus\bZ/2&\dots\\
\bZ/2\oplus\bZ/2\oplus\bZ/2&0&\bZ/2\oplus\bZ/2\oplus\bZ/2 & \dots \\
\bZ/2\oplus\bZ/2&0&\bZ/2\oplus\bZ/2&\dots\\
\bZ/2 & 0 &\bZ/2 &\dots & \imap\\
\bZ/2 & 0 &\bZ/2 &\dots
\end{array}\]
Using the argument of section \ref{sec:finite}, we can show that 
\[Q_{2k,2}([\G^b_0])\in H_{2k}(\cM(S_{2,2};\p);\bZ/2)\]
generates the $\bZ/2$'s in the marked row.
\end{example}

\bibliographystyle{plain}
\bibliography{references}

\begin{thebibliography}{10}

\bibitem{Bodig}
Carl-Friedrich B\"odigheimer.
\newblock Moduli spaces of \up{R}iemann surfaces with boundary.
\newblock To appear, 2003.

\bibitem{Bod_Tillmann}
Carl-Friedrich B{\"o}digheimer and Ulrike Tillmann.
\newblock Stripping and splitting decorated mapping class groups.
\newblock In {\em Cohomological methods in homotopy theory (Bellaterra, 1998)},
  volume 196 of {\em Progr. Math.}, pages 47--57. Birkh\"auser, Basel, 2001.

\bibitem{Bowditch}
B.~H. Bowditch and D.~B.~A. Epstein.
\newblock Natural triangulations associated to a surface.
\newblock {\em Topology}, 27(1):91--117, 1988.

\bibitem{Chas_Sullivan}
Moira Chas and Dennis Sullivan.
\newblock String topology.
\newblock To appear.

\bibitem{Cohen_Tillmann}
F.~R. Cohen and Ulrike Tillmann.
\newblock Toward homology operations for mapping class groups.
\newblock In {\em Homotopy theory via algebraic geometry and group
  representations (Evanston, IL, 1997)}, volume 220 of {\em Contemp. Math.},
  pages 35--46. Amer. Math. Soc., Providence, RI, 1998.

\bibitem{Cohen_Godin}
Ralph Cohen and V\'eronique Godin.
\newblock A polarized view of string topology.
\newblock {\em Proceedings of 2002 Conference on Top. Geo. and Quant. Field
  theory in honor of G. Segal}, 2003.

\bibitem{ConVog}
James Conant and Karen Vogtmann.
\newblock On a theorem of {K}ontsevich.
\newblock {\em Algebr. Geom. Topol.}, 3:1167--1224 (electronic), 2003.

\bibitem{CV86}
Marc Culler and Karen Vogtmann.
\newblock Moduli of graphs and automorphisms of free groups.
\newblock {\em Invent. Math.}, 84(1):91--119, 1986.

\bibitem{ehr}
Ralf Ehrenfried.
\newblock {\em Die {H}omologie der {M}odulr\"aume berandeter {R}iemannscher
  {F}l\"achen von kleinem {G}eschlecht}.
\newblock Bonner Mathematische Schriften [Bonn Mathematical Publications], 306.
  Universit\"at Bonn Mathematisches Institut, Bonn, 1998.
\newblock Dissertation, Rheinische Friedrich-Wilhelms-Universit\"at Bonn, Bonn,
  1997.

\bibitem{Fied_Song}
Zbigniew Fiedorowicz and Yongjin Song.
\newblock The braid structure of mapping class groups.
\newblock {\em Sci. Bull. Josai Univ.}, (Special issue 2):21--29, 1997.
\newblock Surgery and geometric topology (Sakado, 1996).

\bibitem{Galatius_Modp}
S{\o}ren Galatius.
\newblock Mod $p$ homology of the stable mapping class group.
\newblock To appear, 2004.

\bibitem{Harer_Zagier}
J.~Harer and D.~Zagier.
\newblock The {E}uler characteristic of the moduli space of curves.
\newblock {\em Invent. Math.}, 85(3):457--485, 1986.

\bibitem{Harer}
John~L. Harer.
\newblock The virtual cohomological dimension of the mapping class group.
\newblock {\em Inventiones Mathematicae}, 84(1):157--176, 1986.

\bibitem{Harer_Stability}
John~L. Harer.
\newblock Stability of the homology of the moduli spaces of {R}iemann surfaces
  with spin structure.
\newblock {\em Math. Ann.}, 287(2):323--334, 1990.

\bibitem{Igu00}
Kiyoshi Igusa.
\newblock {\em Higher {F}ranz-{R}eidemeister torsion}, volume~31 of {\em AMS/IP
  Studies in Advanced Mathematics}.
\newblock American Mathematical Society, Providence, RI, 2002.
\newblock Chapter 8.

\bibitem{Jones}
John D.~S. Jones.
\newblock Cyclic homology and equivariant homology.
\newblock {\em Invent. Math.}, 87(2):403--423, 1987.

\bibitem{Kauf_Liv_Pen}
Ralph~M. Kaufmann, Muriel Livernet, and R.~C. Penner.
\newblock Arc operads and arc algebras.
\newblock {\em Geom. Topol.}, 7:511--568 (electronic), 2003.

\bibitem{Kontsevich_Intersection}
Maxim Kontsevich.
\newblock Intersection theory on the moduli space of curves and the matrix
  {A}iry function.
\newblock {\em Comm. Math. Phys.}, 147(1):1--23, 1992.

\bibitem{Mad}
Ib~Madsen and Michael Weiss.
\newblock The stable moduli space of \up{R}iemann surfaces : \up{M}umford's
  conjecture.
\newblock preprint, 2002.

\bibitem{Miller}
Edward~Y. Miller.
\newblock The homology of the mapping class group.
\newblock {\em J. Differential Geom.}, 24(1):1--14, 1986.

\bibitem{Milnor}
John Milnor.
\newblock {\em Introduction to algebraic {$K$}-theory}.
\newblock Princeton University Press, Princeton, N.J., 1971.
\newblock Annals of Mathematics Studies, No. 72.

\bibitem{Mumford}
David Mumford.
\newblock Towards an enumerative geometry of the moduli space of curves.
\newblock In {\em Arithmetic and geometry, Vol. II}, volume~36 of {\em Progr.
  Math.}, pages 271--328. Birkh\"auser Boston, Boston, MA, 1983.

\bibitem{Nielsen}
Jakob Nielsen.
\newblock Die isomorphismengruppe der frien gruppen.
\newblock {\em Annals of Mathematics}, 91:169--209, 1924.

\bibitem{Penner}
R.~C. Penner.
\newblock The decorated {T}eichm\"uller space of punctured surfaces.
\newblock {\em Comm. Math. Phys.}, 113(2):299--339, 1987.

\bibitem{Penner_complex}
R.~C. Penner.
\newblock Perturbative series and the moduli space of {R}iemann surfaces.
\newblock {\em J. Differential Geom.}, 27(1):35--53, 1988.

\bibitem{Stasheff}
James~Dillon Stasheff.
\newblock Homotopy associativity of {$H$}-spaces. {I}, {II}.
\newblock {\em Trans. Amer. Math. Soc. 108 (1963), 275-292; ibid.},
  108:293--312, 1963.

\bibitem{Strebel}
Kurt Strebel.
\newblock {\em Quadratic differentials}, volume~5 of {\em Ergebnisse der
  Mathematik und ihrer Grenzgebiete (3) [Results in Mathematics and Related
  Areas (3)]}.
\newblock Springer-Verlag, Berlin, 1984.

\bibitem{Tillmann}
Ulrike Tillmann.
\newblock On the homotopy of the stable mapping class group.
\newblock {\em Invent. Math.}, 130(2):257--275, 1997.

\end{thebibliography}

\end{document}